\documentclass[10pt,oneside]{article}
\input amssymb.sty
\usepackage[cp1250]{inputenc}
\usepackage[T1]{fontenc}
\usepackage{color}
\usepackage{amsmath,amssymb}

\begin{document}
\newcommand{\qed}{\rule{1.5mm}{1.5mm}}
\newcommand{\proof}{\textit{Proof. }}
\newcommand{\ccon}{\rightarrowtail}
\newcommand{\y}{\mathrm{id}}
\newcommand{\C}{\mathbf{C}}
\newtheorem{theorem}{Theorem}[section]
\newtheorem{lemma}[theorem]{Lemma}
\newtheorem{remark}[theorem]{Remark}
\newtheorem{example}[theorem]{Example}
\newtheorem{corollary}[theorem]{Corollary}
\newtheorem{proposition}[theorem]{Proposition}
\newtheorem{claim}[theorem]{Claim}

\begin{center}
{\LARGE\textbf{Approximation of holomorphic maps from Runge domains
to affine algebraic varieties}\vspace*{3mm}}\\
{\large Marcin Bilski\footnote{Department of Mathematics and Computer Science, Jagiellonian University, \L ojasiewicza 6, 30-348 Krak\'ow, Poland; e-mail: marcin.bilski@im.uj.edu.pl} and Adam Parusi\'nski\footnote{Laboratoire J.-A. Dieudonn\'e, Universit\'e de Nice-Sophia Antipolis, Parc Valrose, 06108 Nice Cedex 02, France; e-mail: adam.parusinski@unice.fr

The authors were partially supported by ANR project STAAVF (ANR-2011 BS01 009); 

M.~Bilski was partially supported by the NCN grant 2011/01/B/ST1/03875.}}\vspace*{8mm}\\
\end{center}
\begin{abstract}
\hspace*{-4.9mm}We present a geometric proof of the theorem saying that
holomorphic maps from Runge domains to affine algebraic varieties admit approximation by Nash maps.
Next we generalize this theorem.
\vspace*{1mm}\\
\textbf{MSC (2010):} 32E30, 32S15, 32S45
\end{abstract}
\section{Introduction}
  After the seminal papers of Artin \cite{A68}, \cite{A}
the fundamental problem of algebraic approximation of holomorphic maps satisfying polynomial equations  has been studied by several mathematicians
(see e.g. \cite{BoK2},
\cite{DLS}, \cite{vD}, \cite{Fo}, \cite{Ku}, \cite{Lem}, \cite{TT3}).  The following 
result, which can be viewed as a global version of Artin's approximation theorem, 
is due to L. Lempert (see \cite{Lem}, p. 335).
\begin{theorem}\label{mainapprr}Let $V,W$ be complex affine algebraic varieties,
let $K\subset W$ be a holomorphically convex compact set and let $f:K\rightarrow V$
be a holomorphic map. Then $f$ can be uniformly approximated by a sequence
$f_{\nu}:K\rightarrow V$ of Nash maps.
\end{theorem}
(For the definition of Nash maps see Section \ref{prelnash}; 
In the case of $V$  nonsingular Theorem \ref{mainapprr} had been proved before
in \cite{DLS}. If $W=\mathbf{C}$ and $V$ is arbitrary then it follows from \cite{vD}.)
  Artin's approximation theorem is local and its proof uses 
 Weierstrass Preparation.   The original proof of Lempert's approximation theorem 
 \cite{Lem}, pp 338-339, relies on
the general N\'eron desingularization, 
a deep and difficult result of commutative algebra for
which the reader is referred to \cite{An}, \cite{Og}, \cite{Po1}, \cite{Po2}, \cite{Sp}, 
\cite{Sn}.

Theorem \ref{mainapprr} is expressed in terms of analytic
geometry and has had numerous applications in the theory of several complex variables (see
\cite{B2}, \cite{DF1}, \cite{DF2}, \cite{FT}, \cite{Lem}, \cite{MY}, \cite{NR},
\cite{TT3}).  It is natural to ask whether one can replace N\'eron desingularization
by simpler geometric methods. The main purpose of this paper is to present a new 
geometric proof of Theorem \ref{mainapprr}
based on classical  arguments of singularity theory and complex analysis (see Section \ref{pmt11}).  In the last section, Section \ref{general21}, we show how our method allows us to generalize Lempert's result.

Several variants of Artin's approximation theorem 
turned out to be very useful in singularity theory and complex geometry. It is difficult to give here
a full account of this research.  Instead we refer the reader to two recent papers 
\cite{HR}, \cite{Mi} and references therein.

Our main results will be preceded by an outline of the proof of Theorem~\ref{mainapprr}
where we present the main ideas
and explain why the proof is organized
in the way it is (see Section \ref{sketch}).
The preliminary material is gathered in Section \ref{prelimain}.

\section{Outline of the proof of Theorem \ref{mainapprr}}
\label{sketch} First the problem is reduced (by means of standard methods of multidimensional
complex analysis) to the case where $W=\mathbf{C}^n$ for some integer $n,$ and $K$ is a
compact polydisc (see Section \ref{pmt11}). Then it is sufficient to prove the following
\begin{theorem}\label{maintheorem}
Let $f:U\rightarrow V$ be a holomorphic map,
where $U\subseteq\mathbf{C}^n$ is an open polydisc and $V\subseteq\mathbf{C}^q$ is
an algebraic variety. Then for every open $U_0\Subset U$ there
is a sequence $f_{\nu}:U_0
\rightarrow V$ of Nash maps converging uniformly to $f|_{U_0}.$
\end{theorem}
Theorem \ref{maintheorem} will be obtained by reducing to the case when
$\mathrm{dim}f^{-1}(\mathrm{Sing}(V))<n-1$ and by applying
Lemma \ref{special1}.

Throughout the proof the domain on which the relevant functions
are defined will be shrunk for several times. For this reason in Section \ref{pmt11} we work
with a fixed compact set $K,$ and the functions are defined in some neighborhood of $K$ (which
can be changed). In this outline, for simplicity of notation we assume that $U$ has all the properties
which actually are obtained after shrinking this domain.

Let $f^{-1}(\mathrm{Sing}(V))_{(n-1)}$ denote the union
of all $(n-1)$--dimensional irreducible components of $f^{-1}(\mathrm{Sing}(V))$
and let $\overline{f(U)}^z$ denote the Zariski closure of $f(U)$ (i.e. the smallest algebraic set
containing $f(U)$).
We shall explain the idea of the proof of Theorem~\ref{maintheorem} additionally assuming
that $\overline{f(U)}^z=V,$ and that $f^{-1}(\mathrm{Sing}(V))_{(n-1)}$ has a finite
number of irreducible components. Making
these assumptions we do not lose generality. (This is because first $V$ can be replaced by
$\overline{f(U)}^z.$ Then one can shrink $U$ to obtain
the finiteness condition.) Our aim is to construct a holomorphic map $F_1:U\rightarrow V_1$
such that: $\overline{F_1(U)}^z=V_1$ and $F_1^{-1}(\mathrm{Sing}(V_1))_{(n-1)}$ has
fewer irreducible components than $f^{-1}(\mathrm{Sing}(V))_{(n-1)},$ and if $F_1$ can be
approximated by Nash maps into $V_1$ then $f$ can be approximated by Nash maps into $V.$ When
this is accomplished we can replace $f$ by $F_1$ and repeat the whole process. Such
repetitions lead us to the case solved by Lemma \ref{special1}, application of which finishes
the proof.

Let us start with preliminary remarks. First $V\subset\mathbf{C}^k\times\mathbf{C}^{q-k}$ can
be assumed to be an irreducible normal analytic space with proper projection onto $\mathbf{C}^k,$ where
$k=\mathrm{dim}V.$ A reduction to the case where this assumption is satisfied is standard
(for details see the proof of Theorem \ref{maintheorem} in Section \ref{pmt11}). Then the set
$\Sigma_V\subset\mathbf{C}^k$ (defined in Section \ref{aswpraaa}) is either empty or purely $(k-1)$-dimensional.
(This is because otherwise $V$ would not be locally irreducible contradicting normality.)
Consequently, there is a reduced $N\in\mathbf{C}[w_1,\ldots,w_k]$ such that
$N^{-1}(0)=\Sigma_V.$ Let $G=\pi\circ f,$ where $\pi:\mathbf{C}^k\times\mathbf{C}^{q-k}
\rightarrow\mathbf{C}^k$ denotes the natural projection.
Since $\overline{f(U)}^z=V,$ we have $G(U)\nsubseteq\Sigma_V.$ Therefore
$(N\circ G)^{-1}(0)$ is either purely $(n-1)$--dimensional or empty. We can assume,
shrinking $U$ if needed, that $(N\circ G)^{-1}(0)$ has a finite number of irreducible
components.

Now our main tools are Propositions \ref{redtojm}, \ref{hhhk} proved in Sections
\ref{susect3.1}, \ref{susect3.2}, respectively, and Corollary \ref{maintool}. First,
Proposition \ref{redtojm} enables us to reduce the problem to the case when
\vspace*{2mm}\\
($\flat$)\hspace*{20mm}
${G}^{-1}(\mathrm{Sing}(\Sigma_{V}))_{(n-1)}\subseteq f^{-1}(\mathrm{Sing}(V))_{(n-1)}.$\vspace*{2mm}\\
To be more precise, Proposition \ref{redtojm} provides us with a suitable linear change of the coordinates in
$\mathbf{C}^q$ after which ($\flat$) holds.

Next we will construct a holomorphic map $F_{\ast}$ into an algebraic variety
$V_{\ast}$ with $F_{\ast}^{-1}(\mathrm{Sing}(V_{\ast}))\subseteq {G}^{-1}(\mathrm{Sing}(\Sigma_{V}))$ and
such that if there is a sequence $F_{\ast,\nu}$ of Nash maps into $V_{\ast}$ approximating $F_{\ast},$
then there is a sequence $G_{\nu}$ of Nash maps into $\mathbf{C}^k$
approximating $G$ such that
$\{(N\circ G_{\nu})^{-1}(0)\}$
converges to $(N\circ G)^{-1}(0)$ in the sense of chains.
(We say that a sequence $\{B_{\nu}\}$ of purely $s$--dimensional analytic sets converges
to a purely $s$--dimensional analytic set $B$ in the sense of chains if
$\{Z_{\nu}\}$ converges to $Z,$ where $Z_{\nu}$ and $Z$
are chains obtained by assigning multiplicity $1$ to all irreducible components
of $B_{\nu}$ and $B,$ respectively. For the definition of the convergence of
chains see Section \ref{holchai}.)

Observe that, by
Corollary \ref{maintool}, the existence of $G_{\nu}$ as above implies that
there is a sequence
$f_{\nu}$ of Nash maps into $V$ approximating  $f.$ Moreover, by ($\flat$), the
number of the irreducible components of $F_{\ast}^{-1}(\mathrm{Sing}(V_{\ast}))_{(n-1)}$
does not exceed the number of the irreducible components of $f^{-1}(\mathrm{Sing}(V))_{(n-1)}.$
Since the former number can be equal to the latter one,  in general we cannot
define $V_1=V_{\ast},$ $F_1=F_{\ast}.$ However, $V_{\ast}$ will have a very special description
whose modification will allow us to construct $V_1, F_1$ with all the required properties.

Let us describe how to obtain $V_{\ast}, F_{\ast}$ and $V_1, F_1.$ (Details are presented
in the
proof of Proposition \ref{hhhk}.) Let $A_1,\ldots,A_p$ denote the (pairwise distinct)
irreducible components of $(N\circ{G})^{-1}(0).$ Since $U$ is an open polydisc, we have
$N\circ{G}=u_1^{\alpha_1}\cdot\ldots\cdot u_p^{\alpha_p}\bar{R},$ where
$\bar{R}\in\mathcal{O}(U)$ is a nowhere vanishing function, $u_1,\ldots,u_p$ are minimal
defining functions for $A_1,\ldots,A_p,$ and $\alpha_1,\ldots,\alpha_p$ are positive
integers. (Recall that $u\in\mathcal{O}(U)$ is called a minimal defining function for $A$ if
$A=u^{-1}(0)$ and for every open subset $D\subseteq U$ and $v\in\mathcal{O}(D)$ with $A\cap
D\subseteq v^{-1}(0),$ there is $g\in\mathcal{O}(D)$ such that $v=g\cdot u|_D.$ It is well
known that the existence of minimal defining functions is a consequence of universal
solvability of the second Cousin problem on $U$ which, if
 $U$ is a domain of holomorphy, is equivalent to $H^2(U,\mathbf{Z})=0,$ cf. \cite{Ho}.)

Now define ${F}_{\ast}=({G},u_1,\ldots,u_p,\bar{R}),$
$${V}_{\ast}=\{(w_1,\ldots,w_k,u_1,\ldots,u_p,\bar{R})\in\mathbf{C}^{k+p+1}:
N(w_1,\ldots,w_k)=u_1^{\alpha_1}\cdot\ldots\cdot u_p^{\alpha_p}\bar{R}\},$$
and suppose that there are sequences ${G}_{\nu},u_{1,\nu},\ldots, u_{p,\nu},\bar{R}_{\nu}$
of Nash maps converging locally uniformly to ${G},u_1,\ldots,u_p,\bar{R}$ such that
$N\circ{G}_{\nu}=u_{1,\nu}^{\alpha_1}\cdot\ldots\cdot u_{p,\nu}^{\alpha_p}\bar{R}_{\nu}.$
Since
$u_1,\ldots,u_p$ are minimal defining functions, $\{(N\circ{G}_{\nu})^{-1}(0)\}$
converges to $(N\circ{G})^{-1}(0)$ in the sense of chains. Since $N$ is reduced,
${F}_{\ast}^{-1}(\mathrm{Sing}({V_{\ast}}))
\subseteq{G}^{-1}(\mathrm{Sing}(N^{-1}(0))).$ The functions
$u_1,\ldots,u_p,\bar{R}$ will be chosen in such a way that $\overline{F_{\ast}(U)}^z=V_{\ast}.$

Let us turn to $V_1, F_1.$ If
$F_{\ast}^{-1}(\mathrm{Sing}({V_{\ast}}))_{(n-1)}=\emptyset,$ then
set $V_1={V_{\ast}}, F_1={F_{\ast}}.$ Otherwise one can assume that
$A_1\subseteq F_{\ast}^{-1}(\mathrm{Sing}({V_{\ast}}))_{(n-1)},$ and then
we will construct
$V_1, F_1$ with
$$F_1^{-1}(\mathrm{Sing}(V_1))_{(n-1)}\subseteq
\overline{{G}^{-1}(\mathrm{Sing}(N^{-1}(0)))_{(n-1)}\setminus
A_1}.$$
(For any $B\subseteq\mathbf{C}^{q},$ $\overline{B}$ denotes the closure
in the Euclidean topology.) The construction will be carried out in $\alpha_1$ steps.
More precisely, one step
will be repeated for $\alpha_1$ times, each time with different input data.
In each step we
modify the lefthand
side of the equation $N(w_1,\ldots,w_k)=u_1^{\alpha_1}\cdot\ldots\cdot u_p^{\alpha_p}\bar{R}$
and add to the
system a collection of extra equations of the form $q_j=v_ju_1,$ where $q_j$ are suitably
chosen polynomials
and $v_j$ are new variables.
This operation, which allows us to decrease the power of $u_1$ by $1,$
can be viewed as some sort of blowing-up. After $\alpha_1$ repetitions
we obtain a system of polynomial equations
$N_{\alpha_1}(w_1,\ldots,w_k,v_1,\ldots,v_{t_{\alpha_1}})=
u_2^{\beta_2}\cdot\ldots\cdot u_p^{\beta_p}R_{\alpha_1},$ $q_j=v_ju_1,$ for
$j=1,\ldots,t_{\alpha_1},$ which defines some variety containing $V_1$ as an
irreducible component.

Together with the equations we will introduce new functions corresponding to the variables
$v_j, R_{\alpha_1}$
(also denoted by $v_j, R_{\alpha_1}$) which will become components of the map $F_1.$
\section{Preliminaries}
\label{prelimain}
\subsection{Runge domains and polynomial polyhedra}\label{polyhed}
A domain of holomorphy $\Omega\subset\mathbf{C}^n$ is called a Runge
domain if every function $f\in\mathcal{O}(\Omega)$ can be uniformly
approximated on every compact subset of $\Omega$ by polynomials in
$n$ complex variables.

We say that $P$ is a polynomial polyhedron in $\mathbf{C}^n$ if
there exist polynomials in $n$ complex variables $q_1,\ldots,q_s$
and real constants $c_1,\ldots,c_s$ such that
$$P=\{x\in\mathbf{C}^n:|q_1(x)|\leq c_1,\ldots,|q_s(x)|\leq c_s\}.$$

The following theorem is a straightforward consequence of Theorem
2.7.3 and Lemma 2.7.4 from \cite{Ho}.
\begin{theorem}\label{hor1}Let $\Omega\subset\mathbf{C}^n$ be a Runge
domain. Then for every $\Omega_0\Subset\Omega$ there exists
a compact polynomial polyhedron $P\subseteq\Omega$ such that
$\Omega_0\subseteq P.$
\end{theorem}

The following fact from \cite{Ho} (p. 55) is well known.
\begin{theorem}\label{hor2}Let $f$ be a holomorphic function in a
neighborhood of a compact polynomial polyhedron
$K\subset\mathbf{C}^n.$ Then $f$ can be uniformly approximated on
$K$ by polynomials in $n$ complex variables.
\end{theorem}
\subsection{Nash maps and sets}\label{prelnash}
Let $\Omega$ be an open subset of $\mathbf{C}^n$ and let $f$ be a
holomorphic function on $\Omega.$ We say that $f$ is a Nash
function at $x_0\in\Omega$ if there exist an open neighborhood $U$
of $x_0$ and a polynomial
$P:\mathbf{C}^n\times\mathbf{C}\rightarrow\mathbf{C},$ $P\neq 0,$
such that $P(x,f(x))=0$ for $x\in U.$ A holomorphic function
defined on $\Omega$ is said to be a Nash function if it is a Nash
function at every point of $\Omega.$ A holomorphic mapping defined
on $\Omega$ with values in $\mathbf{C}^N$ is said to be a Nash
mapping if each of its components is a Nash function.

A subset $Y$ of an open set $\Omega\subset\mathbf{C}^n$ is said to
be a Nash subset of $\Omega$ if and only if for every
$y_0\in\Omega$ there exists a neighborhood $U$ of $y_0$ in
$\Omega$ and there exist Nash functions $f_1,\ldots,f_s$ on $U$
such that $$Y\cap U=\{x\in U: f_1(x)=\ldots=f_s(x)=0\}.$$

The following proposition explains the relation between
Nash and algebraic sets (cf. \cite{Tw}).
\begin{proposition}\label{rednash}
Let $X$ be an irreducible Nash subset of an open set
$\Omega\subset\mathbf{C}^n.$ Then there exists an algebraic subset
$Y$ of $\mathbf{C}^n$ such that $X$ is an analytic irreducible
component of $Y\cap\Omega.$ Conversely, every analytic irreducible
component of $Y\cap\Omega$ is an irreducible Nash subset of
$\Omega.$
\end{proposition}
\subsection{Convergence of closed sets and holomorphic chains}\label{holchai}
Let $U$ be an open subset in $\mathbf{C}^m.$ By a holomorphic chain in $U$ we mean a formal
sum $A=\sum_{j\in J}\alpha_jC_j,$ where $\alpha_j\neq 0$ for $j\in J$ are integers and
$\{C_j\}_{j\in J}$ is a locally finite family of pairwise distinct irreducible analytic
subsets of $U$ (see \cite{Ch}, \cite{Tw2}, cf. also \cite{Ba}). The set $\bigcup_{j\in J}C_j$
is called the support of $A$ and is denoted by $|A|$ whereas the sets $C_j$ are called the
components of $A$ with multiplicities $\alpha_j.$ The chain $A$ is called positive if
$\alpha_j>0$ for all $j\in J.$ If all the components of $A$ have the same dimension $n$ then
$A$ will be called an $n-$chain.

Below we introduce the convergence of holomorphic chains in $U$.
To do this we first need the notion of the local uniform
convergence of closed sets. Let $Y,Y_{\nu}$ be closed subsets of
$U$ for $\nu\in\mathbf{N}.$ We say that $\{Y_{\nu}\}$ converges to
$Y$ locally uniformly if:\vspace*{2mm}\\
(1l) for every $a\in Y$ there exists a sequence $\{a_{\nu}\}$ such
that $a_{\nu}\in Y_{\nu}$
and\linebreak\hspace*{7mm}$a_{\nu}\rightarrow a$
in the standard topology of $\mathbf{C}^m,$\\
(2l)  for every compact subset $K$ of $U$ such that $K\cap
Y=\emptyset$ it holds $K\cap Y_{\nu}=\emptyset$\linebreak\hspace*{6.3mm}for almost all $\nu.$\vspace*{1mm}\\
Then we write $Y_{\nu}\rightarrow Y.$ For details concerning the
topology of local uniform convergence see \cite{TwW}.

We say that a sequence $\{Z_{\nu}\}$ of positive $n$-chains
converges to a positive $n$-chain $Z$ if:\vspace*{2mm}\\
(1c) $|Z_{\nu}|\rightarrow |Z|,$\\
(2c) for each regular point $a$ of $|Z|$ and each submanifold $T$
of $U$ of dimension\linebreak\hspace*{7mm}$m-n$ transversal to
$|Z|$ at $a$ such that the closure $\overline{T}$ (in $U$) is compact and
\linebreak\hspace*{7mm}$|Z|\cap\overline{T}=\{a\},$ we have
$deg(Z_{\nu}\cdot T)=deg(Z\cdot T)$ for almost
all $\nu.$\vspace*{2mm}\\
Then we write $Z_{\nu}\ccon Z.$ (By $Z\cdot T$ we denote the
intersection product of $Z$ and $T$ (cf. \cite{Tw2}). (1c) and the choice of $a, T$ in (2c) imply
that the chains $Z_{\nu}\cdot T$ and $Z\cdot T$ for sufficiently large
$\nu$ have finite supports and the degrees are well defined.
Recall that for a chain $A=\sum_{j=1}^d\alpha_j\{a_j\},$
$deg(A)=\sum_{j=1}^d\alpha_j.$)

When we say that a sequence $\{X_{\nu}\}$ of purely
$n$--dimensional analytic sets converges to a purely $n$-dimensional
analytic set $X$ \textit{in the sense of chains}, we mean that
the sequence $\{Z_{\nu}\}$ of $n$--chains converges to the $n$--chain $Z,$
where $Z_{\nu}, Z$ are obtained by assigning the multiplicity $1$ to all
irreducible components of $X_{\nu}, X$ respectively.

\subsection{Analytic sets with proper projection}\label{aswpraaa}
Let $\pi:\mathbf{C}^m\times\mathbf{C}^s\rightarrow\mathbf{C}^m$ be
the natural projection, let $\Omega$ be a domain in
$\mathbf{C}^m,$ and let $Y$ be a purely $m$--dimensional analytic
subset of $\Omega\times\mathbf{C}^s$ such that
$\pi|_Y:Y\rightarrow\Omega$ is a proper map.
By $\Sigma_Y$ we denote the set of all points
$a\in\Omega$ such that the fiber of $\pi|_Y$ over $a$ does not have the maximal cardinality.
Recall that $\Sigma_Y$ (called the discriminant of $\pi|_Y$) is an analytic subsets of
$\Omega$ (cf. \cite{Ch}).

For algebraic sets, we need a slightly more general notion.
Let $\mathcal{L}(\mathbf{C}^N,\mathbf{C}^m)$ denote the vector space of all linear maps from
$\C^N$ to $\C^m.$ Let $V\subset\C^N$ be algebraic of pure dimension $m$ and let
$A\in\mathcal{L}(\C^N,\C^m)$ such that $A|_{V}:V\rightarrow\C^m$ is a proper map. By
$\Sigma_A\subset\C^m$ we denote the set of points $a\in\C^m$
such that the fiber of $A|_V$ over $a$ does not have the maximal cardinality. Recall that
$\Sigma_A$ (called the discriminant of $A|_V$) is algebraic.
Set
$s_A:=(A|_V)^{-1}(\mathrm{Sing}(\Sigma_A)).$
When $A$ is the natural projection from $\C^{N}=\C^{m}\times\C^{s}$ to $\C^m,$ we often write
$\Sigma_V$ instead of $\Sigma_A.$
\subsection{Holomorphic maps into algebraic varieties}
For any subset $B$ of $\mathbf{C}^q$ let $\overline{B}^{z}$ denote the Zariski
closure of $B$ i.e. the intersection of all algebraic subvarieties
of $\mathbf{C}^q$ containing $B.$ For any algebraic subvariety $V$ of $\mathbf{C}^q$
by $I(V)$ we denote the ideal of all polynomials $p\in\mathbf{C}[y_1,\ldots,y_q]$
such that $V\subseteq p^{-1}(0).$ For any $g_1,\ldots,g_s\in\mathbf{C}[y_1,\ldots,y_q]$
by $I(g_1,\ldots,g_s)$ we denote the ideal generated by $g_1,\ldots, g_s.$

\begin{lemma}\label{monter}Let $\Omega$ and $B$ be a domain in $\mathbf{C}^n$
and an irreducible analytic subset of $\Omega,$ respectively, and let
$g=(g_1,g_2)\in\mathcal{O}(\Omega,\mathbf{C}^q\times\mathbf{C}^r).$
Let
$\delta,h_1,\ldots,h_{t_1}\in\mathbf{C}[y_1,\ldots,y_q]$ be such that
$$\overline{g_1(B)}^z\setminus\delta^{-1}(0)=\{y\in\mathbf{C}^q\setminus
\delta^{-1}(0):h_j(y)=0\mbox{ for }j=1,\ldots,t_1\},$$
where $t_1=q-\mathrm{dim}\overline{g_1(B)}^z,$ and
$\delta|_{\overline{g_1(B)}^z}\neq 0$
and for every $a\in\overline{g_1(B)}^z\setminus\delta^{-1}(0),$
the map $(h_1,\ldots,h_{t_1}):\mathbf{C}^q\rightarrow\mathbf{C}^{t_1}$
is a submersion in some neighborhood of $a$ in $\mathbf{C}^q,$ and
$\delta I(\overline{g_1(B)}^z)\subseteq I(h_1,\ldots,h_{t_1}).$
Then there are
$t_2-t_1$ polynomials $h_{t_1+1},\ldots,h_{t_2}\in\mathbf{C}
[y_1,\ldots,y_q,y_{q+1},\ldots,y_{q+r}]$, where $t_2=q+r-\mathrm{dim}\overline{g(B)}^z,$
and there is
$\hat{\delta}\in\mathbf{C}[y_1,\ldots,y_q,y_{q+1},\ldots,y_{q+r}]$
such that
$$\overline{g(B)}^z\setminus\hat{\delta}^{-1}(0)=\{y\in\mathbf{C}^q\times\mathbf{C}^r\setminus
\hat{\delta}^{-1}(0):h_j(y)=0\mbox{ for }j=1,\ldots,t_2\},$$
and $\hat{\delta}|_{\overline{g(B)}^z}\neq 0,$ and
for every $b\in\overline{g(B)}^z\setminus\hat{\delta}^{-1}(0),$
the map $(h_1,\ldots,h_{t_2}):\mathbf{C}^{q+r}\rightarrow\mathbf{C}^{t_2}$
is a submersion in some neighborhood of $b$ in $\mathbf{C}^{q+r},$
and $\hat{\delta} I(\overline{g(B)}^z)\subseteq I(h_1,\ldots,h_{t_2}).$
\end{lemma}
\textit{Proof of Lemma \ref{monter}.} Let us denote $C_1=\overline{g_1(B)}^z,$
$C_2=\overline{g(B)}^z.$ Since $B$
is irreducible, $C_2$ is irreducible as well.
Then there are
$\hat{\delta}_1, \hat{h}_{1},\ldots,\hat{h}_{t_2}\in\mathbf{C}[y_1,\ldots,y_{q+r}],$
such that
$$C_2\setminus\hat{\delta}_1^{-1}(0)=\{y\in\mathbf{C}^q\times\mathbf{C}^r\setminus
\hat{\delta}_1^{-1}(0):\hat{h}_j(y)=0\mbox{ for }j=1,\ldots,t_2\},$$
and $\hat{\delta}_1|_{C_2}\neq 0,$ and
for every $b\in C_2\setminus\hat{\delta}_1^{-1}(0),$
the map $(\hat{h}_1,\ldots,\hat{h}_{t_2}):\mathbf{C}^{q+r}\rightarrow\mathbf{C}^{t_2}$
is a submersion in some neighborhood of $b$ in $\mathbf{C}^{q+r},$
and for every
${G}\in\mathbf{C}[y_1,\ldots,y_{q+r}]$ with $C_2\subseteq{G}^{-1}(0)$
there are $\hat{r}_1,\ldots,\hat{r}_{t_2}\in\mathbf{C}[y_1,\ldots,y_{q+r}]$
such that $\hat{\delta}_1\cdot{G}=\sum_{j=1}^{t_2}\hat{r}_j\hat{h}_j.$
(See \cite{Lo}, pp. 402-405.)

Let us show that there are ${h}_{t_1+1},\ldots,{h}_{t_2}\in\{\hat{h}_1,\ldots,\hat{h}_{t_2}\}$
with the required properties.
Observe that $C_2\subseteq C_1\times\mathbf{C}^r,$
which implies that $\hat{\delta}_1\cdot h_{i}=\sum_{j=1}^{t_2}b_{j,i}\hat{h}_{j},$
for $i=1,\ldots,t_1,$ where $b_{j,i}\in\mathbf{C}[y_1,\ldots,y_{q+r}].$ Next,
$(\hat{\delta}_1\circ g)\cdot(\delta\circ g_1)|_B\neq 0.$
Indeed, otherwise either $\hat{\delta}_1|_{C_2}=0$
or $\delta|_{C_1}=0.$ Consequently, there is $x_0\in B$
such that $(\hat{\delta}_1h_1,\ldots,\hat{\delta}_1h_{t_1})$ is
a submersion in a neighborhood of $g(x_0),$ and therefore there
are $j_1,\ldots,j_{t_1}$ such that the determinant $d(y_1,\ldots,y_{q+r})$
of the matrix
$[b_{j_k,i}(y_1,\ldots,y_{q+r})]_{k=1,\ldots,t_1;i=1,\ldots,t_1}$
satisfies $d(g(x_0))\neq 0.$ This implies that $d|_{C_2}\neq 0$
and there are $c_{k,i}, d_{l,i}\in\mathbf{C}[y_1,\ldots,y_{q+r}],$
for $k,i=1,\ldots,t_1$ and $l\in J=\{1,\ldots,t_2\}\setminus\{j_1,\ldots,j_{t_1}\}$
such that
$d\cdot\hat{h}_{j_i}=\sum_{k=1}^{t_1}c_{k,i}h_k+\sum_{l\in J}d_{l,i}\hat{h}_l.$

Now it is clear that
the assertion of the lemma is satisfied with $\hat{\delta}=d\cdot\hat{\delta}_1\cdot\delta$
and
$(h_1,\ldots,h_{t_1},h_{t_1+1},\ldots,h_{t_2}),$ where
$h_k=\hat{h}_{j_{k}},$ for $k=t_1+1,\ldots,t_2$ where
$\{j_{t_1+1},\ldots,j_{t_2}\}=J.$
\qed
\begin{remark}\label{1222}{\em Let $a_1,\ldots,a_p\in\mathbf{C}^{q+r}\setminus C_2.$
Then\vspace*{1mm}\\
(a) In the first paragraph of the proof of Lemma \ref{monter},
$\hat{h}_1,\ldots,\hat{h}_{t_2},$ $\hat{\delta}_1$ can be chosen so that
$\hat{h}_i(a_j)\neq 0$ $\neq\hat{\delta}_1(a_j)$ for $i=1,\ldots,t_2$ and $j=1,\ldots,p.$\\
(b) If $h_{i}(a_j)\neq 0\neq\delta(a_j)$ for $i=1,\ldots,t_1$ and
$j=1,\ldots,p$ then $h_{t_1+1},\ldots,h_{t_2}$ and the $b_{j,i}$'s can be chosen so that
$h_i(a_j)\neq 0\neq\hat{\delta}(a_j)$ for $i=1,\ldots,t_2$ and $j=1,\ldots,p.$}
\end{remark}
\textit{Proof of Remark \ref{1222}.} As for (a), it is sufficient to prove the assertion
with $C_2, a_1,\ldots,a_p$ replaced by $\Phi(C_2),\Phi(a_1),\ldots,\Phi(a_p),$ where $\Phi$
is any linear isomorphism. Therefore we may assume, applying a linear change of the variables
in $\C^{q+r}$ if needed, that $\pi|_{C_2}:C_2\rightarrow\C^{\mathrm{dim}(C_2)+1}$ is a proper
map, where $\pi$ denotes the projection onto the first $\mathrm{dim}(C_2)+1$ coordinates
of $\C^{q+r}.$ Moreover, $\pi(a_j)\notin\pi(C_2)$ for every $j=1,\ldots,p$  and
the fiber of $\pi|_{C_2}$ over $a$ consists of one element for every
$a\in\mathrm{Reg}(\pi(C_2)).$ We may also assume that $\rho|_{\pi(C_2)}:\pi(C_2)\rightarrow
\C^{\mathrm{dim}(C_2)}$ is a proper map, where $\rho$ denotes the projection onto the first
$\mathrm{dim}(C_2)$ coordinates of $\C^{\mathrm{dim}(C_2)+1}$ and that the fibers of
$\rho|_{\pi(C_2)}$ have maximal cardinality over $\rho(\pi(a_j))$ for every $j=1,\ldots,p.$

Let $P\in(\C[y_1,\ldots,y_{\mathrm{dim}(C_2)}])[y_{\mathrm{dim}(C_2)+1}]$ be the irreducible
monic polynomial with $P^{-1}(0)=\pi(C_2).$ Then it is well known
(cf. \cite{Lo}, pp. 402-405) that one can take $\hat{\delta}_1,\hat{h}_1,\ldots,\hat{h}_{t_2}$ such that
$\hat{h}_1=P$ and $\hat{\delta}_1$ is (a power of) the discriminant of $P.$
Of course, $\hat{h}_1,\hat{\delta}_1$ satisfy the requirements. Finally, if needed,
we can replace $\hat{h}_j,$ for $j\geq 2,$ by $\hat{h}_j+\epsilon_j\hat{h}_1$ where
$\epsilon_j\in\C$ ($|\epsilon_j|$ small) to obtain (a).

Let us turn to (b) (assuming that we have (a)). In view of (a), the fact that
$h_{j},$ for $j=t_1+1,\ldots,t_2,$
are chosen among $\hat{h}_1,\ldots,\hat{h}_{t_2}$, and that
$\hat{\delta}=d\cdot\hat{\delta}_1\cdot\delta,$ it is sufficient to observe that
the $b_{j,i}$'s can be chosen in such a way that the determinant $d(a_{j})\neq 0$ for $j=1,\ldots,p.$
This however is obvious because for every $j\notin J,$ $b_{j,i}$ can be replaced by
$b_{j,i}+\epsilon_{j,i}\hat{h}_l$ for any $l\in J$ and $\epsilon_{j,i}\in\C$ with $|\epsilon_{j,i}|$
small.\qed
\begin{lemma}\label{modbaseq}Let $K$ be a compact polydisc in $\mathbf{C}^n$ and let
$G\in\mathcal{O}(K,\mathbf{C}^k)$ such that $\overline{G(K)}^z=\mathbf{C}^k.$
Let $0\neq N\in\mathbf{C}[y_1,\ldots,y_k]$ and $u_1,\ldots,u_p, R\in\mathcal{O}(K)$ satisfy
$N\circ G=u_1^{\alpha_1}\cdot\ldots\cdot u_p^{\alpha_p}R$ for some positive integers
$\alpha_1,\ldots,\alpha_p.$
Then there are nowhere vanishing functions
$v_1,\ldots,v_p\in\mathcal{O}(K)$ and there is $S\in\mathcal{O}(K)$ such that
$N\circ G=(u_1v_1)^{\alpha_1}\cdot\ldots\cdot(u_pv_p)^{\alpha_p}S$ and
$\overline{(G,u_1v_1,\ldots,u_pv_p)(K)}^z=\mathbf{C}^{k+p}.$
\end{lemma}
\proof
It is sufficient to show that there are nowhere vanishing
$v_1,\ldots,v_p$ $\in\mathcal{O}(K)$ such that
$\overline{(G,u_1v_1,\ldots,u_pv_p)(K)}^z=\mathbf{C}^{k+p},$ by which the other asser\-tions
follow immediately. In other words, it is sufficient to show that if
$g_1,\ldots,g_t\in\mathcal{O}(K)$ are algebraically independent over $\mathbf{C},$ 
then for any $u\in\mathcal{O}(K),$ $u\neq 0,$
there is a nowhere vanishing $v\in\mathcal{O}(K)$ such that $g_1,\ldots,g_t,u\cdot v$ are
also algebraically independent.

We have two cases: either $g_1,\ldots,g_t, u$ are algebraically independent (and there
is nothing to prove) or not. In the latter case it is sufficient to
show that there is a nowhere vanishing $v\in\mathcal{O}(K)$ such that $g_1,\ldots,g_t,v$
are algebraically independent over $\mathbf{C}$ because then $g_1,\ldots,g_t,u\cdot v$ are also such
(cf. \cite{Lan} for basic facts on algebraic extensions).
Define a family of one-variable functions: $\chi_1(x)=\exp(x)$ and $\chi_j(x)=\chi_{j-1}(\exp(x))$ for $j>1.$
Then $\{\chi_j\}_{j=1}^{k}$ are algebraically independent for any $k$.
Hence, if $g_1,\ldots,g_t$ are algebraically independent,
there is $j$ such that $g_1,\ldots,g_t,\chi_j$ are also independent
(where $\chi_j$ is treated as an $n$-variable function; cf. \cite{Lan} for basic facts on transcendental extensions).\qed
\subsection{A discriminant criterion for the existence of algebraic approximations}
\label{adcfeaa} Let us recall a result from \cite{B4} which is one of the main
tools in the present paper. Let $U\subset\mathbf{C}^n$ be a domain and let
$\pi:U\times\mathbf{C}^k\rightarrow U$ denote the natural projection. Let $X\subset
U\times\mathbf{C}^k$ be an analytic subset of pure dimension $n$ with proper projection onto
$U.$ Recall that $s(\pi|_{X})$ denotes the cardinality of the generic fiber in $X$ over $U.$
For any analytic set $C,$ $C_{(n-1)}$ denotes the union of all $(n-1)$--dimensional
irreducible components of $C.$ For the definition of $\Sigma_X$ see Section \ref{aswpraaa}.
\begin{theorem}\label{component} Let $\{X_{\nu}\}$ be a sequence
of purely $n$-dimensional analytic subsets of
$U\times\mathbf{C}^k$ with proper projection onto $U$ converging
locally uniformly to $X$ such that $s(\pi|_{X})=s(\pi|_{X_{\nu}})$
for $\nu\in\mathbf{N}.$ Assume that
$\{(\Sigma_{X_{\nu}})_{(n-1)}\}$ converges to
$(\Sigma_X)_{(n-1)}$ in the sense of holomorphic
chains. Then for every analytic subset $Y$ of
$U\times\mathbf{C}^k$ of pure dimension $n$ such that $Y\subseteq X$
and for every open relatively compact subset $\tilde{U}$ of $U$
there exists a sequence $\{Y_{\nu}\}$ of purely $n$-dimensional
analytic subsets of $\tilde{U}\times\mathbf{C}^k$ converging to
$Y\cap(\tilde{U}\times\mathbf{C}^k)$ in the sense of holomorphic
chains such that $Y_{\nu}\subseteq X_{\nu}$ for every
$\nu\in\mathbf{N}.$
\end{theorem}
The assumption that $\{(\Sigma_{X_{\nu}})_{(n-1)}\}$ converges to $(\Sigma_X)_{(n-1)}$ in the
sense of chains is essential. It is not difficult to observe that otherwise $X_{\nu}$ can be
(and usually are) irreducible even when $X$ is not. Then, in particular, $X$ can contain the
graph $Y$ of some map holomorphic on $\tilde{U}$ whereas $X_{\nu}$ do not contain any such
graphs. To prove Theorem \ref{maintheorem}, we will use Theorem \ref{component} in the case when
$Y$ is the graph of a holomorphic map to be approximated by Nash maps.

Let $V$ be a purely $m$--dimensional algebraic variety in $\mathbf{C}^m\times\mathbf{C}^k$
with proper projection onto $\mathbf{C}^m.$ Assume that $\Sigma_V=N^{-1}(0),$ where $N$ is a
polynomial in $m$ variables. Let $f:U\rightarrow V$ be a holomorphic map such that
$\tilde{f}(U)\nsubseteq N^{-1}(0)$ where $\tilde{f}$ is the map consisting of the first $m$
components of $f.$
Theorem \ref{component} implies the following
\begin{corollary}\label{maintool}
Let $\tilde{f}_{\nu}\in\mathcal{O}(U,\mathbf{C}^m)$ be a sequence of Nash maps
converging to $\tilde{f}$ locally uniformly such that $\{(N\circ\tilde{f}_{\nu})^{-1}(0)\}$
converges to $(N\circ\tilde{f})^{-1}(0)$ in the sense of chains. Then for every analytic
subset $Y$ of $U\times\mathbf{C}^k$ of pure dimension $n$ such that
$Y\subseteq (\tilde{f}\times\mathrm{id_{\C^k}})^{-1}(V)$
and for every open relatively compact subset $\tilde{U}$ of $U$ there is a sequence
$\{Y_{\nu}\}$ of purely $n$--dimensional Nash subsets of $\tilde{U}\times\mathbf{C}^k$
converging to $Y\cap(\tilde{U}\times\mathbf{C}^k)$ in the sense of holomorphic chains
such that $Y_{\nu}\subseteq(\tilde{f}_{\nu}|_{\tilde{U}}\times\mathrm{id_{\C^k}})^{-1}(V)$
for every $\nu\in\mathbf{N}.$
In particular, there is a sequence $f_{\nu}\in\mathcal{O}(\tilde{U},V)$ of Nash maps
converging to $f|_{\tilde{U}}$ uniformly.
\end{corollary}
\proof Only the last sentence requires explanation. Let $\hat{f}$ be the map consisting
of the last $k$ components of $f.$ Fix open $\tilde{U}\Subset\hat{U}\Subset U.$ Since
$\mathrm{graph}(\hat{f})\subseteq
(\tilde{f}\times\mathrm{id_{\C^k}})^{-1}(V),$ there are
purely $n$--dimensional Nash sets $Y_{\nu}\subseteq(\tilde{f}_{\nu}|_{\hat{U}}\times\mathrm{id_{\C^k}})^{-1}(V)$
such that $\{Y_{\nu}\}$ converges to $\mathrm{graph}(\hat{f})\cap(\hat{U}\times\C^k)$
in the sense of chains.

The fibers of
$(\tilde{f}_{\nu}|_{\hat{U}}\times\mathrm{id_{\C^k}})^{-1}(V)\subseteq\hat{U}\times\mathbf{C}^k$ over $\hat{U}$
are uniformly bounded by some bound independent of $\nu$
(which follows by the fact that $V$ has proper projection onto $\C^m$ and $\hat{U}\Subset U$).
Hence the fibers of $Y_{\nu}$ over $\hat{U}$ are also such. This
implies, in view of the fact that $\{Y_{\nu}\}$ converges in the
sense of chains to the graph of a map, that the fibers of $Y_{\nu}\cap(\tilde{U}\times\C^k)$
are singletons for almost all $\nu.$
In other words, $Y_{\nu}\cap(\tilde{U}\times\C^k)$ is the graph of some (Nash) map
$\hat{f}_{\nu}:\tilde{U}\rightarrow\C^k$ for almost all $\nu.$ Now we can define
${f}_{\nu}=(\tilde{f}_{\nu}|_{\tilde{U}},\hat{f}_{\nu}).$\qed

\section{Proof of Theorem \ref{mainapprr}}\label{pmt11}
\subsection{Reduction to the case when $K$ is a polydisc}
A proof of Theorem \ref{mainapprr} can be reduced (by means of standard techniques of
multidimensional complex analysis) to the case where $W=\mathbf{C}^n,$ $K$ is a
compact polynomial
polyhedron, and $f:D\rightarrow V,$ where $D$ is an open neighborhood of $K$ in
$\mathbf{C}^n$ (see \cite{Lem}, p. 339). Let us assume that this reduction has been done.
As we will show below, a modification of the method
presented in \cite{Lem} (p. 339) allows us to assume that $K$ is a compact polydisc.

First, for some $m,$ there are a compact polydisc $E$ in $\mathbf{C}^{n+m},$ polynomials
${P}_1,\ldots,{P}_m$ and a mapping $F$ of a neighborhood $U$ of $E$ into $\mathbf{C}^q$ such
that in some neighborhood $\tilde{D}$ of $K$ in $D$ we have $f(z)=F(z,P_1(z),\ldots,P_m(z))$
(cf. \cite{Ho}, p. 55). Moreover, if $Z\subset\mathbf{C}^n\times\mathbf{C}^m$ denotes the
graph of the map $z\mapsto(P_1(z),\ldots,P_m(z)),$ then $\tilde{K}=(K\times\mathbf{C}^m)\cap
Z\subseteq E$ is a polynomially convex compact subset of $Z.$ We can take $U$ to be an open
polydisc. It is clear that in order to approximate $f$ it is sufficient to approximate the
map $F|_{\tilde{K}}:\tilde{K}\rightarrow V$ by Nash maps into $V.$

Let $Q_1,\ldots,Q_r$ be polynomials in $q$ complex variables such that
$V=\{Q_1=\ldots=Q_r=0\}.$ Let $\tilde{Z}$ denote the union of all analytic irreducible
components of $Z\cap U$ which have a non-empty intersection with
$\tilde{D}\times\mathbf{C}^m.$ Observe that, for $i=1,\ldots,r,$ $Q_i\circ F$ vanishes
identically on $\tilde{Z}.$ Pick $\alpha\in\mathcal{O}(U)$ non-vanishing identically on any
analytic irreducible component of $\tilde{Z}$ but vanishing identically on the other
irreducible components of $Z\cap U$ (cf. \cite{Ho}, p. 192). Then $\alpha(y)\cdot
Q_i(F(y))=0$ for every $y\in Z\cap U$ and $i=1,\ldots,r.$ Therefore, in view of $Z$ being
algebraic, there are polynomials $R_1,\ldots, R_t$ in $n+m$
complex variables, vanishing identically on $Z\cap U,$ and such that\vspace*{2mm}\\
($\natural$)\hspace*{10mm}$\alpha(y)\cdot Q_i(F(y))=\sum_{j=1}^tb_{i,j}(y)R_j(y), \mbox{
for }y\in U\mbox{ and }i=1,\ldots,r$\vspace*{2mm}\\
with certain holomorphic functions
$b_{i,j}.$

In the space $\mathbf{C}^{1+q+n+m+rt}$ with coordinates $x_1,
w_1,\ldots,w_q,u_1,\ldots,u_{n+m},$ $v_{1,1},v_{1,2},\ldots,v_{r,t}$ consider the variety $T$
defined by the equations
$$x_1\cdot Q_i(w)=\sum_{j=1}^tv_{i,j}R_j(u), \mbox{ for }i=1,\ldots,r.$$
By ($\natural$), the image of the map $$g: U\ni y\mapsto
(\alpha(y),F(y),y,b_{i,j}(y))\in\mathbf{C}^{1+q+n+m+rt}$$ is contained in $T.$

Now suppose
that there is an open polydisc $U'$ with $E\subset U'\subset U$ such that
$g$ can be approximated on $U'$ by a Nash map
$g'(y)=(\alpha'(y), F'(y), y'(y), b_{i,j}'(y))$ whose image is contained in $T.$
If this approximation is close enough, then $y'(y)$ has the
inverse $\tilde{y}$ on $E$ close to the identity. Replacing $g'(y)$ by $g'(\tilde{y}(y)),$ we
can assume that $y'(y)=y.$
Consequently, $\alpha'(y)\cdot
Q_i(F'(y))=0$ for every $y\in\tilde{Z}\cap U'$ and every $i.$ But $\alpha'(y)$ does not vanish
identically on any irreducible component of $\tilde{Z}\cap U'$ if the approximation is close enough.
Therefore $Q_i(F'(y))=0$ for every $y\in\tilde{Z}\cap U'$ and every $i,$ which implies that
$F'(\tilde{K})\subset V.$

Thus in order to obtain the required Nash approximation $F'|_{\tilde{K}}$ of $F|_{\tilde{K}}$
it suffices to approximate $g:E\rightarrow T,$ where $E$ is a compact polydisc. Now, we see
that to prove Theorem \ref{mainapprr}, it is sufficient to prove Theorem \ref{maintheorem}.
\subsection{Proof of Theorem \ref{maintheorem}}
First we will focus on the case when
$\mathrm{dim}f^{-1}(\mathrm{Sing}(V))\leq n-2.$

\begin{lemma}\label{special1}
Let $f:U\rightarrow V$ be a holomorphic map,
where $U\subseteq\mathbf{C}^n$ is a Runge domain and $V\subseteq\mathbf{C}^q$ is
an algebraic variety. Assume that $\mathrm{dim}f^{-1}({\mathrm{Sing(V)}})<n-1.$
Then for every open $U_0\Subset U$ there is a sequence $f_{\nu}:U_0
\rightarrow V$ of Nash maps converging uniformly to $f|_{U_0}.$
\end{lemma}
\proof For an elementary proof of the lemma for $n=1$ the reader is referred to
\cite{vD}. Let us assume that $n\geq 2.$
For any $\C$-linear subspace $L$ of $\C^n$ let $L^{\bot}$ denote the orthogonal
complement of $L$ in $\C^n.$ Fix an open set $U_0\Subset U.$

Since $\mathrm{dim}f^{-1}(\mathrm{Sing}(V))<n-1,$ there are $\epsilon>0,$ $(n-2)$--dimensional
linear subspaces
$L_1,\ldots,L_t\subset\C^n$ and open bounded balls $B_j\subset L_j$ and $B_j'\subset L_j^{\bot},$ for
$j=1,\ldots,t,$ such that $P_j:=B_j+B_j'\Subset U$ and $\overline{U}_0\subseteq\bigcup_{j=1}^tP_j,$
and
$(\overline{B_{j}}+\overline{B'_{j,\epsilon}})\cap f^{-1}(\mathrm{Sing}(V))= \emptyset,$
for $j=1,\ldots,t,$ where $B'_{j,\epsilon}=\{x\in B'_j:\mathrm{dist}(x,\partial
B'_j)<\epsilon\}.$

Observe that for every $i\neq j$ such that $P_i\cap P_j\neq\emptyset$ there are open balls
${P}_{i,j}\subseteq B_i,$ $P_{j,i}\subseteq B_j$ and open connected sets $l_{i,j}\subseteq
B_i',$ $l_{j,i}\subseteq B_j'$ such that $\overline{P_{i,j}+ l_{i,j}}\cap
f^{-1}(\mathrm{Sing}(V))=\emptyset,$ $\overline{P_{j,i}+ l_{j,i}}\cap
f^{-1}(\mathrm{Sing}(V))=\emptyset,$ $l_{i,j}\cap B'_{i,\epsilon}\neq\emptyset,$ $l_{j,i}\cap
B'_{j,\epsilon}\neq\emptyset,$ and $(P_{i,j}+l_{i,j})\cap(P_{j,i}+l_{j,i})\neq\emptyset.$
Indeed, pick $z\in(P_i\cap P_j) \setminus f^{-1}(\mathrm{Sing}(V)).$ We have
$z=u_i+v_i=u_j+v_j$ for some $u_i\in L_i,$ $v_i\in L_i^{\bot},$ $u_j\in L_j,$ $v_j\in
L_j^{\bot}.$ Next pick $v_i'\in B_{i,\epsilon}',$ $v_j'\in B_{j,\epsilon}'.$ Since
$u_i+v_i', u_j+v_j'\notin f^{-1}(\mathrm{Sing}(V)),$ there are a path $k_{i,j}\subset B_i'$
connecting $v_i,v_i'$ and a path $k_{j,i}\subset B_j'$ connecting $v_j,v_j'$ such that
$(u_i+k_{i,j})\cap f^{-1}(\mathrm{Sing}(V))=\emptyset=(u_j+k_{j,i})\cap f^{-1}(\mathrm{Sing}(V)).$
Now it suffices to take $P_{i,j}, P_{j,i}$ to be small balls centered at
$u_i, u_j$ in $B_i, B_j,$ respectively, and $l_{i,j},$ $l_{j,i}$ to be small neighborhoods
of $k_{i,j},$ $k_{j,i}$ in $B_i',B_j',$ respectively.

Define $E=\bigcup_{j=1}^t({B_j}+ B'_{j,\epsilon}) \cup\bigcup_{i,j=1}^t(P_{i,j}+ l_{i,j}),$
assuming that $P_{i,j}+ l_{i,j}=\emptyset$ if $P_i\cap P_j=\emptyset.$ By the facts that $U$
is a Runge domain, $E\Subset U$ and $f(\overline{E})\cap\mathrm{Sing}(V)=\emptyset,$ there is
a sequence ${f}_{\nu}:E\rightarrow V$ of Nash maps approximating $f|_E$ uniformly (cf.
\cite{DLS} p. 334; the idea of the proof is as follows:
since
$\overline{f(E)}\cap\mathrm{Sing}(V)=\emptyset,$ there are an open neighborhood $N$ of
$\overline{f(E)}$ in $\C^q$ together with a Nash retraction $\tau$ of $N$ onto some open
neighborhood of $\overline{f(E)}$ in $V.$ Now $f|_E$ can be
approximated by polynomial maps into $\C^q$ whose
restrictions to $E$ have images in $N.$ These restrictions can be composed with
$\tau$ yielding the required Nash approximations of $f|_E$).

Observe that if $f_{\nu}$ has
a holomorphic extension to $\bigcup_{j=1}^tP_j,$
then the proof will be completed. Indeed, by the maximum principle,
if such $f_{\nu}$ approximates $f$ on $E$ then it also approximates $f$
on $\bigcup_{j=1}^tP_j.$
Moreover, if $f_{\nu}$ is a holomorphic map on $\bigcup_{j=1}^tP_j$ and
a Nash map on $E$ then it is a Nash map on $\bigcup_{j=1}^tP_j.$

For every $j$ put $E_j=(B_j+B'_{j,\epsilon})\cup\bigcup_{k=1}^t(P_{j,k}+l_{j,k})$
(again assuming $P_{j,k}+l_{j,k}=\emptyset$ if $P_j\cap P_k=\emptyset$).
The Hartogs extension theorem implies that
for every $j,$ $f_{\nu}|_{E_j}$ has an
extension $f_{j,\nu}:P_j\rightarrow V$ such that
$f_{j,\nu}|_{\{z\}+B'_{j}}$
is a holomorphic map for every $z\in B_j.$
But then, since $f_{j,\nu}|_{B_j+B'_{j,\epsilon}}$
is a holomorphic map, the Cauchy integral formula implies
that $f_{j,\nu}$ is a continuous separately holomorphic map.
Hence it is a holomorphic map.

It remains to show that for every $i,j,$
$f_{i,\nu}|_{P_i\cap P_j}=f_{j,\nu}|_{P_i\cap P_j}.$
Fix $i\neq j$ such that $P_i\cap P_j\neq\emptyset.$
Since
$C=(P_{i,j}+l_{i,j})\cap(P_{j,i}+l_{j,i})\subseteq E_i\cap E_j\subseteq P_i\cap P_j$
we have $f_{j,\nu}|_{C}=f_{\nu}|_{C}=f_{i,\nu}|_{C}.$
But $C\neq\emptyset$ and $P_i\cap P_j$ is connected
so $f_{i,\nu}|_{P_i\cap P_j}=f_{j,\nu}|_{P_i\cap P_j}$
and the proof is complete.\qed\vspace*{2mm}\\
\textbf{Notation.}
Let $K$ be a connected compact subset of $\mathbf{C}^n$
such that
$\mathrm{int}K\neq\emptyset$ and let
$g:K\rightarrow C\subseteq\mathbf{C}^q.$ By
$g_D:D\rightarrow C$ we denote a holomorphic map
such that $g_D|_K=g,$ where $D$ is an open connected neighborhood of $K$.
If $g_D$ exists, then $g$ is called holomorphic. If $g_{D}$ is a Nash map then $g$ is
called a Nash map. The collection of all holomorphic maps from $K$
to $C$ will be denoted by $\mathcal{O}(K,C).$ For $C=\mathbf{C}$
we write $\mathcal{O}(K).$ A sequence $g_{\nu}\in\mathcal{O}(K,C),$
for $\nu\in\mathbf{N},$ is said to converge to $g\in\mathcal{O}(K,C)$
uniformly if there is an open $D'\supset K$ for which there are
$g_{D'}, g_{\nu,D'},$ $\nu\in\mathbf{N},$
such that $g_{\nu,D'}$ converges to $g_{D'}$ uniformly.
Let $h\in\mathcal{O}(D,\mathbf{C}^q)$ for some open $D\subset\C^n.$
Let $Y\subset\mathbf{C}^q$ be an analytic set. Then by
$h^{-1}(Y)_{(n-1)}$ we denote the union of all $(n-1)$--dimensional irreducible components of
$h^{-1}(Y).$\vspace*{2mm}\\
\textit{Proof of Theorem \ref{maintheorem}.} Fix an open $U_0\Subset U$
(which clearly can be assumed to be connected) and a compact polydisc $K$ with $U_0\subset K\subset U.$
One can assume that $\overline{f(K)}^z=V$ (because otherwise
$V$ can be replaced by $\overline{f(K)}^z$), and that
$(f|_D)^{-1}(\mathrm{Sing}(V))_{(n-1)}\neq\emptyset$ for every open neighborhood $D$ of $K$
(as otherwise Lemma \ref{special1} finishes the proof).

Put $F_0=f|_{{K}}, V_0=V.$ We iterate the following process starting from $F_0.$ Suppose we
have $F_i\in\mathcal{O}({K},V_i)$ such that
$\overline{F_i(K)}^z=V_i,$
$F_{i,D}^{-1}(\mathrm{Sing}(V_i))_{(n-1)}$ $\neq\emptyset,$ for every open neighborhood
$D$ of $K,$
where $V_i\subset\mathbf{C}^{q_i}.$ We will show that
there is $F_{i+1}\in\mathcal{O}({K}, V_{i+1}),$
where $V_{i+1}\subset\mathbf{C}^{q_{i+1}}$ is an algebraic variety, such that:\vspace*{2mm}\\
(x) if there is a sequence $F_{i+1,\nu}\in\mathcal{O}(K,V_{i+1})$ of Nash maps
converging uniformly to $F_{i+1},$ then there is a sequence
$F_{i,\nu}\in\mathcal{O}(K,V_i)$ of Nash maps converging uniformly
to $F_i,$\\
(y) $\overline{F_{i+1}(K)}^z=V_{i+1}$ and there are an open $D\supset K$ and an
irreducible component $T$ of $(F_{i,D})^{-1}(\mathrm{Sing}(V_{i}))_{(n-1)}$
with $T\cap K\neq\emptyset$
such that
$$(F_{i+1,D})^{-1}(\mathrm{Sing}(V_{i+1}))_{(n-1)}\subseteq
\overline{(F_{i,D})^{-1}(\mathrm{Sing}(V_{i}))_{(n-1)}\setminus T}.$$

Let us show that once $F_{i+1}$ is constructed, the proof will be completed.
Set $C_{i,D'}:=(F_{i,D'})^{-1}(\mathrm{Sing}(V_{i}))_{(n-1)}.$
Let $I_i\subset\mathcal{O}(K)$ be the ideal of all
$\alpha\in\mathcal{O}(K)$ such that $\alpha_{D'}|_{C_{i,D'}}=0$ for some open
$D'\supset K.$ It is well known that for every analytic hypersurface $H$
in an open polydisc $D'$ there is $g\in\mathcal{O}(D')$ such that $H=g^{-1}(0)$ (cf. \cite{Ho}). This
fact and (y) imply that $I_{i}\varsubsetneq I_{i+1}.$ Therefore if
$C_{i,D'}\neq\emptyset$ for every $i$ and
every open $D'\supset K,$ then there is an infinite ascending sequence of ideals in
$\mathcal{O}(K).$ But $\mathcal{O}(K)$ is noetherian (cf. \cite{Fr}) so
there must be $i_0$ and an open $D'\supset K$ such that $C_{i_0,D'}=\emptyset$
(i.e. $\mathrm{dim}F_{i_0,D'}^{-1}(\mathrm{Sing}(V_{i_0}))<n-1$).
Now it is clear that Lemma \ref{special1} allows us to complete the proof
if, given $F_i,$ we can construct $F_{i+1}$ satisfying (x) and (y).

Put $k_i=\mathrm{dim}(V_i).$ Let us show how to
construct $F_{i+1}.$ First observe that
$V_i$ is irreducible (because $K$ is a polydisc and $V_i=\overline{F_i(K)}^z$). 
We can also assume that $V_i$ is a normal analytic space.
Indeed, if $V_i$ is not normal then we can replace $V_i, F_i$ by
$\tilde{V}_i,\tilde{F}_i,$ where $\pi:\tilde{V}_i\rightarrow V_i$ is the normalization of
$V_i,$ whereas $\tilde{F}_i:{K}\rightarrow\tilde{V}_i$ is a holomorphic map such that
$\pi\circ\tilde{F}_i=F_i.$ (The existence of $\tilde{F}_i$ is an immediate consequence of the
fact that
$\pi|_{\tilde{V}_i\setminus\pi^{-1}(\mathrm{Sing}(V_i))}:\tilde{V}_i\setminus\pi^{-1}(\mathrm{Sing}(V_i))
\rightarrow V_i\setminus\mathrm{Sing}(V_i)$ is a biholomorphism (see \cite{Lo}, pp 343-346).)

After this preparation let us construct $F_{i+1}.$ Our main tools are Propositions
\ref{redtojm} and \ref{hhhk} (whose proofs are postponed to Sections \ref{susect3.1},
\ref{susect3.2}, respectively) and Corollary \ref{maintool}. For definitions of $\Sigma_A,$ $s_A$ and
$\mathcal{L}(\mathbf{C}^N,\mathbf{C}^m)$ see Section \ref{aswpraaa}.

\begin{proposition}\label{redtojm}
Let $V$ be an algebraic subset of $\C^N$ of pure dimension $m,$ and let
$U\subset\C^n$ be an open polydisc.
Let $f: U\rightarrow V$ be a holomorphic map such that $\overline{f(U)}^z=V.$ Then for every
open $U_0\Subset U$ there is $A\in\mathcal{L}(\C^N,\C^m)$ such that
$A|_{V}:V\rightarrow\mathbf{C}^m$ is a proper map and
$\mathrm{dim}(f|_{U_0})^{-1}(s_A\setminus\mathrm{Sing}(V))\leq n-2. $
\end{proposition}
By Proposition \ref{redtojm}, there is a linear
$A:\C^{q_i}\rightarrow\C^{k_i}$ such that $A|_{V_i}:V_i\rightarrow\C^{k_i}$ is proper
and\vspace*{2mm}\\
(a) $(A|_{V_i}\circ F_{i,D})^{-1}(\mathrm{Sing}(\Sigma_{A}))_{(n-1)}
\subseteq F_{i,D}^{-1}(\mathrm{Sing}(V_i))_{(n-1)},$ for every sufficiently small
open $D\supset K.$\vspace*{2mm}\\
Let $\Phi:\C^{q_i}\rightarrow\C^{q_i}$ be a linear automorphism
such that $A=\pi\circ\Phi,$ where $\pi:\mathbf{C}^{k_i}\times\mathbf{C}^{q_i-k_i}
\rightarrow\mathbf{C}^{k_i}$ is the natural projection. Then $\Sigma_A=\Sigma_{\Phi(V_i)}.$
Since $\Phi(V_i)$ is a normal space, $\Sigma_{\Phi(V_i)}$ is purely
$(k_i-1)$--dimensional or empty.
Therefore there is $N\in\mathbf{C}[w_1,\ldots,w_{k_i}]$ such that
$N^{-1}(0)=\Sigma_{\Phi(V_i)}.$
Set ${G}=\pi\circ\Phi\circ F_i.$ Then (a) can be rewritten as follows\vspace*{2mm}\\
(0) ${G}_D^{-1}(\mathrm{Sing}(N^{-1}(0)))_{(n-1)}\subseteq F_{i,D}^{-1}
(\mathrm{Sing}(V_i))_{(n-1)},$ for every sufficiently small
open $D\supset K.$\vspace*{2mm}

On the other hand, the facts that $\mathrm{dim}(\Phi(V_i))=k_i,$ $\Phi(V_i)$ has proper projection onto $\C^{k_i}$ and
$\overline{\Phi(F_i(K))}^z=\Phi(V_i)$ imply that $\overline{G(K)}^z=\C^{k_i}.$
Hence, ${G}, N$ satisfy the hypotheses of the following
\begin{proposition}\label{hhhk}
Let $E\subset\mathbf{C}^n$ be a compact
polydisc with $\mathrm{int}(E)\neq\emptyset.$
Let $G\in\mathcal{O}(E, \mathbf{C}^k),$ $N\in\mathbf{C}[w_1,\ldots,w_k]$ satisfy
$\overline{G(E)}^z=\C^k,$ $N\neq 0.$
Then
there are an algebraic subset $\tilde{V}$ of some $\mathbf{C}^q$ and
$\tilde{f}\in\mathcal{O}(E,\tilde{V})$ with $\overline{\tilde{f}(E)}^z=\tilde{V}$
such that:\vspace*{2mm}\\
$(1)$ either $\tilde{f}_D^{-1}(\mathrm{Sing}(\tilde{V}))_{(n-1)}=\emptyset$ for some
open $D\supset E$
or there are an open $D\supset E$ and an irreducible component $T$ of
$G_D^{-1}(\mathrm{Sing}(N^{-1}(0)))_{(n-1)}$ with $T\cap E\neq\emptyset$ such that
 $\tilde{f}_D^{-1}(\mathrm{Sing}(\tilde{V}))_{(n-1)}\subseteq\overline{
 G_D^{-1}(\mathrm{Sing}(N^{-1}(0)))_{(n-1)}\setminus T},$\\
$(2)$ if there is a sequence
$\tilde{f}_{\nu}\in\mathcal{O}({E},\tilde{V})$ of Nash maps
converging uniformly to $\tilde{f},$ then there are a sequence
$G_{\nu}\in\mathcal{O}({E},\mathbf{C}^k)$
of Nash maps converging uniformly to $G$ and an open $D'\supset E$
such that $\{(N\circ G_{\nu,D'})^{-1}(0)\}$ converges to $(N\circ G_{D'})^{-1}(0)$ in the sense
of chains.
\end{proposition}

By Proposition \ref{hhhk}, there are an algebraic subset ${V}_{i+1}$ of some
$\mathbf{C}^{q_{i+1}}$ and $F_{i+1}\in\mathcal{O}({K},{V_{i+1}})$ with
$\overline{F_{i+1}(K)}^z=V_{i+1}$
such that:\vspace*{2mm}\\
(3) either
$F_{i+1,D}^{-1}(\mathrm{Sing}({V_{i+1}}))_{(n-1)}=\emptyset$ for some open $D\supset K$
or there are an open $D\supset K$ and an irreducible
component $T$ of ${G}_D^{-1}(\mathrm{Sing}(N^{-1}(0)))_{(n-1)}$ with $T\cap K\neq\emptyset$
such that
$F_{i+1,D}^{-1}(\mathrm{Sing}({V_{i+1}}))_{(n-1)}\subseteq\overline{
{G}_D^{-1}(\mathrm{Sing}(N^{-1}(0)))_{(n-1)}\setminus T},$\vspace*{2mm}\\
(4) if there is a sequence
${F}_{i+1,\nu}\in\mathcal{O}(K,{V}_{i+1})$ of Nash maps
converging uniformly to $F_{i+1},$
then there are a sequence ${G}_{\nu}\in\mathcal{O}(K,\mathbf{C}^{k_i})$
of Nash maps converging uniformly to ${G}$ and an open $D'\supset K$
such that $\{(N\circ{G}_{\nu,D'})^{-1}(0)\}$ converges to
$(N\circ{G}_{D'})^{-1}(0)$ in the sense of chains.\vspace*{2mm}\\
\hspace*{5mm}Now, by (4) and Corollary \ref{maintool}, if there is a sequence
${F}_{i+1,\nu}\in\mathcal{O}(K,{V}_{i+1})$ of Nash maps converging uniformly to $F_{i+1},$
then there is a sequence $\bar{F}_{\nu}\in\mathcal{O}(K,\Phi(V_i))$ of Nash maps converging uniformly to
$\Phi\circ F_i,$ which clearly implies that (x) is satisfied. As for (y), it is an immediate
consequence
of (0) and (3). Thus the proof is complete.\qed\vspace*{2mm}
\subsection{Proof of Proposition \ref{redtojm}}\label{susect3.1}

We follow the notation introduced in subsection 3.4.  
Throughout the proof we fix a nonempty Zariski open subset
$T$ of $\mathcal{L}(\C^N,\C^m)$ such that 
$$\pi:V\times T\rightarrow\C^m\times T, \quad \pi(x,A)=(A(x),A)$$ 
is proper and 
$(\Sigma_{\pi})\cap(\C^m\times\{A\})=\Sigma_A\times\{A\},$ for all $A\in T,$ where
$\Sigma_{\pi}$ denotes the discriminant of
$\pi$. Then by  Bertini Theorem (see for instance \cite{Hartshorne77}
 Corollary 10.9 and Remark 10.9.2, pp. 274-275) replacing $T$ by a smaller nonempty Zariski open subset of $\mathcal{L}(\C^N,\C^m)$,
 if necessary, we have 
$$\mathrm{Sing}(\Sigma_{\pi})\cap(\C^m\times\{A\})=\mathrm{Sing}(\Sigma_A)\times\{A\}, 
\quad \text{ for all }
A\in T.$$
Therefore, if we denote $s_\pi = \pi^{-1} (\mathrm{Sing}(\Sigma_{\pi}))$, then\vspace*{2mm}\\
(a)\hspace*{20mm}$s_\pi \cap(V\times\{A\})= s_A\times\{A\}, 
\quad \text{ for all }
A\in T.$\vspace*{2mm}\\
Since $\mathrm{dim Sing}(\Sigma_A)\leq m-2$ and $A|_V:V\rightarrow\C^m$ is proper
for $A\in T,$ we have\vspace*{2mm}\\
(b)\hspace*{25mm}$\mathrm{dim}(s_A)\leq m-2,$ for all $A\in T.$\vspace*{4mm}\\
For a line $L$ in $\mathcal{L}(\C^N,\C^m)$ we put 
$$
s_L:=\overline{\bigcup_{A\in L\cap T}s_A}^z. 
$$
We claim that $\dim s_L\le m-1$.  Indeed, $\dim (s_\pi \cap (V\times(L\cap T)) \le m-1$ by (a) and (b).  
The image in $V$ of the standard projection 
$V\times T\to V$ of $s_\pi \cap (V\times(L\cap T))$ 
is algebraically constructible (cf. e.g. \cite{Lo}, p. 395) and $s_L$ is its Zariski closure.  This shows the claim 
(cf. e.g. \cite{Lo}, pp. 393-394).  
Finally, for each $x\in\mathrm{Reg}(V)$ there is an $A\in T$ such that 
 $A(x)\notin\Sigma_A$, and hence the set of such $A$ is Zariski open dense, so $\bigcap_{A\in T}s_A\subset\mathrm{Sing}(V)$.  
 Since 
the Zariski topology is noetherian 
\vspace*{2mm}\\
(c)\hspace*{3mm}there are $k\in\mathbf{N}$ and $A_1,\ldots,A_k\in T$ such that
$s_{A_1}\cap\ldots\cap s_{A_k}\subset \mathrm{Sing}(V).$\vspace*{2mm}\\

Now fix an open polydisc $U_0\Subset U$ and set
$Y_A:=(f|_{U_0})^{-1}(s_A\setminus\mathrm{Sing}(V))$.   
For $i\ge 1$ consider the following statement\vspace*{2mm}\\ 
($\mathrm{c}_i$)\hspace*{15mm}$\text{there are } A_1,\ldots ,A_i\in T \text{ such that } 
\mathrm{dim}(\bigcap_{j=1}^i Y_{A_j})\leq n-2.$\vspace*{2mm}\\
By (c),  ($\mathrm{c}_k$) holds.  We will prove that $(\mathrm{c}_i) \Rightarrow   (\mathrm{c}_{i-1})$ for $i\ge 2$,
thus showing ($\mathrm{c}_1$) and hence Proposition 4.2.  

Thus suppose that there are $A_1,\ldots ,A_i\in T \text{ such that } \mathrm{dim}(\bigcap_{j=1}^i Y_{A_j})\leq n-2$.  Let $L$ be the line in $\mathcal{L}(\C^N,\C^m)$ 
containing $A_{i-1}$ and $A_{i}$ and let   
$$Y_L:=(f|_{U_0})^{-1}(s_L\setminus\mathrm{Sing}(V)).$$ 
Since $\overline{f(U)}^z=V$ and $\mathrm{dim}s_L\leq m-1$ we 
obtain\vspace*{2mm}\\
(d)\hspace*{35mm}$\mathrm{dim}(Y_L) \le \dim (f|_{U_0})^{-1}(s_L) \leq n-1.$\vspace*{2mm}

Let $\mathcal{Z}$ denote the finite family of all $(n-1)$--dimensional analytic irreducible
components of $Y_L$. 
For each $Z\in\mathcal{Z}$ the set $\mathcal{A}_Z$ of such $A\in L\cap T$ that\vspace*{2mm}\\
(e)\hspace*{45mm}$Z\subset\bigcap_{j=1}^{i-2}Y_{A_j}\cap Y_A$\vspace*{2mm}\\
is Zariski closed, 
hence either finite or equal $L\cap T$. Indeed, by definitions of $Y_A$ and $Y_L,$
$\mathcal{A}_Z$ equals the set of such $A\in L\cap T$ that $\overline{f(Z)}^z\subset\bigcap_{j=1}^{i-2} s_{A_j}\cap s_A.$ But, by (a), $\bigcap_{j=1}^{i-2} s_{A_j}\cap s_A$ depends algebraically on $A.$

If $\mathcal{A}_Z$ is finite for every $Z\in\mathcal{Z}$ 
then there is $A\in L\cap T$ such that (e) fails for every 
$Z\in\mathcal{Z}$ and then $\dim \bigcap_{j=1}^{i-2}Y_{A_j}\cap Y_A\le n-2$ 
that completes the proof. Thus suppose that there is $Z\in\mathcal{Z}$ for which (e) 
holds for every $A\in 
L\cap T$.  Then 
$$Z\subset\bigcap_{j=1}^{i-2}Y_{A_j}\cap Y_{A_{i-1}} \cap Y_{A_i}$$
that contradicts the assumption $\dim (\bigcap_{j=1}^i Y_{A_j})\leq n-2$.  
\qed


\subsection{Proof of Proposition \ref{hhhk}}\label{susect3.2}

\begin{remark}\label{single1}\emph{
The letters $v_i,$ $u_i,$ $R_i,$ $\hat{w}_i,$ used below denote
either (tuples of) varia\-bles or (tuples of) functions in $x.$
It will be clear from the context whether a given letter denotes a variable
or a function. When we write that a tuple of functions satisfies
some equation, we mean that the equation holds true if every variable
is replaced by the function denoted by the same letter.
}
\end{remark}
If $(N\circ G)^{-1}(0)=\emptyset$ then define
$\tilde{V}\subset\mathbf{C}^{k+1}$ by the equation $N(w_1,\ldots,w_k)$ $=R_0$ and
take  $\tilde{f}(x)=(G(x),R_0(x))=
(G(x),N(G(x))).$ Clearly, $\tilde{V},\tilde{f}$ satisfy the requirements.

Let us assume that $(N\circ G)^{-1}(0)\neq\emptyset.$ Clearly it is sufficient to prove the
proposition in the case where $N$ is reduced which we also assume.

There are an open polydisc $D\supset E$ and a holomorphic
extension $G_D$ of $G.$ Let $A_1,\ldots,A_p$ be all irreducible components of $(N\circ
G_D)^{-1}(0)$ intersecting $E$ and let
$u_1,\ldots,u_p\in\mathcal{O}(D)$ be minimal defining functions for $A_1,\ldots,A_p$
respectively. (Recall that $u\in\mathcal{O}(D)$ is called a minimal defining function for $A$
if $A=u^{-1}(0)$ and for every open subset $U\subseteq D$ and $v\in\mathcal{O}(U)$ with $A\cap
U\subseteq v^{-1}(0),$ there is $g\in\mathcal{O}(U)$ such that $v=g\cdot u|_U.$ It is well
known that the existence of minimal defining functions is a consequence of universal
solvability of the second Cousin problem on $D$ which, if $D$ is a domain of holomorphy, is
equivalent to $H^2(D,\mathbf{Z})=0,$ cf. \cite{Ho}.) Then there are $R_0\in\mathcal{O}(D)$
and positive integers $k_1,\ldots,k_p$ such that $N(G_{D}(x))=u_1(x)^{k_1}\cdot\ldots\cdot
u_p(x)^{k_p}R_0(x)$ and $A_{l}\nsubseteq R_{0}^{-1}(0),$ for $l=1,\ldots,p.$ By Lemma
\ref{modbaseq}, $u_1,\ldots,u_p, R_0$ can be chosen in such a way that
$\overline{(G,u_1,\ldots,u_p)(E)}^z=\C^{k+p}.$

Set $Z_s=G_{D}^{-1}(\mathrm{Sing}(N^{-1}(0)))_{(n-1)},$
$Z_r=\overline{G_{D}^{-1}(N^{-1}(0))\setminus Z_s}.$ If $A_j\nsubseteq Z_s$ for every $j,$
then, after shrinking $D$ if needed, we obtain $Z_s=\emptyset.$ Then define
$\tilde{V}\subset\mathbf{C}^{k+p+1}$ by the equation
$N(w_1,\ldots,w_k)=u_1^{k_1}\cdot\ldots\cdot u_p^{k_p}R_0$ and take $\tilde{f}(x)=(
G_{D}(x),u_1(x),\ldots,u_p(x),R_0(x)).$ Observe that $\tilde{V}, \tilde{f}$ satisfy (2) of
Proposition \ref{hhhk} because $u_1,\ldots,u_p$ are minimal defining functions and
$R_0^{-1}(0)\cap E=\emptyset.$ As for (1), we will show that
$\tilde{f}^{-1}(\mathrm{Sing}(\tilde{V}))_{(n-1)}=\emptyset.$
Suppose that $\tilde{f}(C)\subset\mathrm{Sing}(\tilde{V})$
for some $(n-1)$--dimensional analytic $C\subset D.$ Then
$N\circ G_D|_{C}=0=\frac{\partial N}{\partial w_j}\circ G_D|_{C},$ for $j=1,\ldots,k.$
But $N$ is reduced so $G_D(C)\subset \mathrm{Sing}(N^{-1}(0)).$ This contradicts the fact that $Z_s=\emptyset.$

If $A_j\subseteq Z_s$ for some $j,$ then we
may assume, renumbering the components, that $j=1.$
Put $\hat{w}_1=(w_1,\ldots,w_k),$ $\hat{w}_1(x)=G_{D}(x).$
Now the construction, the aim of which is to remove $A_1$ from $Z_s,$ consists of $k_1$ steps.
\vspace*{2mm}\\
\textit{Step 1.}
Define $C_1=\overline{\hat{w}_1(A_{1})}^z.$ Then $C_1\varsubsetneq\mathbf{C}^k$
is irreducible because
$C_1\subseteq N^{-1}(0)$ and $A_{1}$ is irreducible.
By Lemma \ref{monter} and \cite{Lo}, pp. 402-405, there are $\delta_1,
q_{1},\ldots,q_{t_1}\in\mathbf{C}[\hat{w}_1],$ where $t_1=k-\mathrm{dim}C_1,$ such that
$$C_1\setminus\delta_1^{-1}(0)=
\{\hat{w}_1\in\mathbf{C}^k\setminus\delta_1^{-1}(0):q_{1}(\hat{w}_1)=\ldots=q_{t_1}(\hat{w}_1)=0\},$$
and
$\delta_1|_{C_1}\neq 0,$ and for every $a\in C_1\setminus\delta_1^{-1}(0)$ the map
$(q_{1},\ldots,q_{t_1}):\mathbf{C}^k\rightarrow\mathbf{C}^{t_1}$ is a submersion in some
neighborhood of $a$ in $\mathbf{C}^k.$ Moreover,
$\delta_1 I(C_1)\subseteq I(q_1,\ldots,q_{t_1}).$
Every irreducible component $Z$ of $\overline{\bigcup_{j=1}^p(u_j^{-1}(0))\setminus Z_s}$
satisfies
$\hat{w}_1(Z)\nsubseteq C_1$ because $C_1\subseteq\mathrm{Sing}(N^{-1}(0))$ and
$\hat{w}_1(Z)\nsubseteq\mathrm{Sing}(N^{-1}(0)).$ Therefore, in view of Remark \ref{1222},
we may assume that
every such component $Z$ satisfies
$\hat{w}_1(Z)\nsubseteq\delta_1^{-1}(0).$

The inclusion $C_1\subseteq N^{-1}(0)$ implies that $\delta_{1}N=\sum_{j=1}^{t_1}q_{j}r_{1,j},$
where $r_{1,j}\in\mathbf{C}[\hat{w}_1]$ and
the fact that $u_1$ is a minimal defining function implies that there is
$v_{j}\in\mathcal{O}({D})$ such that
$q_{j}(\hat{w}_1(x))$ $=v_{j}(x)u_1(x)$
for $j=1,\ldots,t_1.$

Let $\hat{v}_1$ denote the tuple $(v_{1},\ldots,v_{t_1})$
of $t_1$ variables. Define
$N_1\in\mathbf{C}[\hat{w}_1,\hat{v}_1]$ by the formula
$$N_1(\hat{w}_1,\hat{v}_1)=\sum_{j=1}^{t_1}v_{j}r_{1,j}(\hat{w}_1)$$
and observe that
$$N_1(\hat{w}_1(x),\hat{v}_1(x))= u_1(x)^{k_1-1}u_2(x)^{k_{2,2}}\cdot\ldots\cdot
u_p(x)^{k_{p,2}}R_{1}(x),$$
where $R_{1}\in\mathcal{O}(D)$ satisfies
$A_{l}\nsubseteq R_1^{-1}(0)$ for $l=1,\ldots,p.$

Let
$V_1\subset\mathbf{C}^{k+t_1+p+1}$
be the algebraic variety defined by the system of equations (in the variables
$\hat{w}_1,\hat{v}_1,u_1,\ldots,u_p,R_1$):\vspace*{2mm}\\
(E,$1$)\hspace*{29.3mm}$N_1(\hat{w}_1,\hat{v}_1)=u_1^{k_{1}-1}u_2^{k_{2,2}}\ldots u_p^{k_{p,2}}R_{1},$\\
(F,$j$)\hspace*{30mm}$q_{j}(\hat{w}_1)=v_{j}u_1, \mbox{ }$ for $j=1,\ldots,t_1.$ \vspace*{2mm}\\
Put $\hat{w}_2(x)=(\hat{w}_1(x),\hat{v}_1(x)),$ and define $g_1\in\mathcal{O}(D)$
by
$$g_1(x)=(\hat{w}_2(x),u_1(x),\ldots,u_p(x),R_1(x)).$$
If $k_{1}=1$
then $\tilde{f}_{D}=g_1, \tilde{V}=\overline{g_1(E)}^z\subseteq V_1$ satisfy the requirements
(see Claims \ref{property1}, \ref{property2}). Otherwise we go to Step 2.\vspace*{2mm}\\
\textit{Step 2.}
Define $C_2=\overline{\hat{w}_2(A_{1})}^z.$ Then $C_2\varsubsetneq\mathbf{C}^{k+t_1}$
is irreducible because
$C_2\subseteq N_1^{-1}(0)$ and $A_{1}$ is irreducible.
Then by Lemma \ref{monter} there are $\delta_2,
q_{t_1+1},\ldots,q_{t_2}\in\mathbf{C}[\hat{w}_2],$ where $t_2=k+t_1-dim C_2,$ such
that
$$C_2\setminus\delta_2^{-1}(0)=\{\hat{w}_2\in\mathbf{C}^{k+t_1}\setminus\delta_2^{-1}(0):
q_{1}(\hat{w}_2)=\ldots=q_{t_2}(\hat{w}_2)=0\},$$
and
$\delta_2|_{C_2}\neq 0,$ and for every $a\in C_2\setminus\delta_2^{-1}(0)$ the map
$(q_{1},\ldots,q_{t_2}):\mathbf{C}^{k+t_1}\rightarrow\mathbf{C}^{t_2}$ is a submersion in some
neighborhood of $a$ in $\mathbf{C}^{k+t_1}.$ Moreover,
$\delta_2 I(C_2)\subseteq I(q_1,\ldots,q_{t_2}).$
Every irreducible component $Z$ of $\overline{\bigcup_{j=1}^p(u_j^{-1}(0))\setminus Z_s}$
satisfies
$\hat{w}_2(Z)\nsubseteq C_2$ because $C_2\subseteq C_1\times\C^{t_1}$ and $\hat{w}_1(Z)\nsubseteq C_1.$
Therefore, in view of Remark \ref{1222}, we may assume that
every such component $Z$ satisfies
$\hat{w}_2(Z)\nsubseteq\delta_2^{-1}(0).$

The inclusion $C_2\subseteq N_1^{-1}(0)$ implies $\delta_2N_1=\sum_{j=1}^{t_2}q_jr_{2,j},$
where $r_{2,j}\in\mathbf{C}[\hat{w}_2],$
and the fact that $u_1$ is a minimal defining function implies that there is
$v_{j}\in\mathcal{O}({D})$ such that
$q_{j}(\hat{w}_2(x))$ $=v_{j}(x)u_1(x)$ for $j=t_1+1,\ldots,t_2.$

Let $\hat{v}_2$ denote the tuple
$(v_{t_1+1},\ldots,v_{t_2})$
of $t_2-t_1$ variables. Define
$N_2\in\mathbf{C}[\hat{w}_2,\hat{v}_2]$ by the formula
$$N_2(\hat{w}_2,\hat{v}_2)=\sum_{j=1}^{t_2}v_{j}r_{2,j}(\hat{w}_2)$$
and observe that
$$N_2(\hat{w}_2(x),\hat{v}_2(x))= u_1(x)^{k_1-2}u_2(x)^{k_{2,3}}\cdot\ldots\cdot
u_p(x)^{k_{p,3}}R_{2}(x),$$
where $R_{2}\in\mathcal{O}(D)$ satisfies
$A_{l}\nsubseteq R_2^{-1}(0)$
for $l=1,\ldots,p.$

Let
$V_2\subset\mathbf{C}^{k+t_2+p+1}$
be the algebraic variety defined by the system of equations
(in the variables $\hat{w}_2,\hat{v}_2,u_1,\ldots,u_p,R_2$):\vspace*{2mm}\\
(E,$2$)\hspace*{28.9mm}$N_2(\hat{w}_2,\hat{v}_2)=
u_1^{k_{1}-2}u_2^{k_{2,3}}\ldots u_p^{k_{p,3}}R_{2},$\\
(F,$j$)\hspace*{30mm}$q_{j}(\hat{w}_2)
=v_{j}u_1, \mbox{ }$ for $j=1,\ldots,t_2.$ \vspace*{2mm}\\
(For $j=1,\ldots,t_1,$ the polynomial $q_j$ is precisely the one from
Step 1 and it does not really depend on $\hat{v}_1.$)
Put $\hat{w}_3(x)=(\hat{w}_2(x),\hat{v}_2(x))$ and define $g_2\in\mathcal{O}(D)$
by
$$g_2(x)=(\hat{w}_3(x),u_1(x),\ldots,u_p(x),R_2(x)).$$
If $k_{1}=2$
then $\tilde{f}_{D}=g_2, \tilde{V}=\overline{g_2(E)}^z\subseteq V_2$ satisfy the requirements
(see Claims \ref{property1}, \ref{property2}). Otherwise we go to Step 3.\vspace*{2mm}\\
Let us describe Step i+1,
assuming that $k_1>i$ and we have completed Step i (i$\geq 2$) after which there are an algebraic
subvariety $V_i\subseteq\mathbf{C}^{k+t_i+p+1}$ and a holomorphic
map $g_i:D\rightarrow V_i$
such that the following hold:\vspace*{2mm}\\
$\bullet$ $$g_i(x)=(\hat{w}_{i+1}(x),u_1(x),\ldots,u_p(x),R_i(x)),$$
$\hat{w}_{i+1}(x)=(\hat{w}_i(x),\hat{v}_{i}(x))\in\mathbf{C}^{k+t_i},$ and
$\hat{v}_i(x)=(v_{t_{i-1}+1}(x),\ldots,$ $v_{t_i}(x)).$ Moreover, $A_{l}\nsubseteq
R^{-1}_i(0)$ for $l=1,\dots,p.$
\vspace*{2mm}\\
$\bullet$ $V_i$ is defined by the equations (in the variables $\hat{w}_i,\hat{v}_i,
u_1,\ldots,u_p,R_i$):\vspace*{2mm}\\
(E,$i$)\hspace*{29.8mm}$N_i(\hat{w}_i,\hat{v}_i)=
u_1^{k_{1}-i}u_2^{k_{2,i+1}}\ldots u_p^{k_{p,i+1}}R_{i},$\\
(F,$j$)\hspace*{30mm}$q_{j}(\hat{w}_i)=v_{j}u_1,\mbox{ }$
for $j=1,\ldots,t_i,$\vspace*{2mm}\\
where $0\neq N_i\in\mathbf{C}[\hat{w}_i,\hat{v}_i],$ and $q_j\in\mathbf{C}[\hat{w}_i]$
for $j=1,\ldots, t_i.$\vspace*{2mm}\\
$\bullet$
There is $\delta_i$ $\in\mathbf{C}[\hat{w}_i]$ such that for
$C_i=\overline{\hat{w}_i(A_{1})}^z\varsubsetneq\mathbf{C}^{k+t_{i-1}}$ the following
hold:
$$C_i\setminus\delta_i^{-1}(0)=\{\hat{w}_i
\in\mathbf{C}^{k+t_{i-1}}\setminus\delta_i^{-1}(0):
q_{j}(\hat{w}_i)=0\mbox{ for }j=1,\ldots,t_i\},$$
and
$\delta_i|_{C_i}\neq 0,$ and for every $a\in C_i\setminus\delta_i^{-1}(0)$ the map
$(q_{1},\ldots,q_{t_i}):\mathbf{C}^{k+t_{i-1}}\rightarrow\mathbf{C}^{t_i}$ is a submersion in some
neighborhood of $a$ in $\mathbf{C}^{k+t_{i-1}},$ and
$\delta_i I(C_i)\subseteq I(q_1,\ldots,q_{t_i}).$
Furthermore, $\hat{w}_{i}(Z)\nsubseteq C_i$ and
$\hat{w}_{i}(Z)\nsubseteq\delta_i^{-1}(0)$
for every irreducible component $Z$ of
$\overline{\bigcup_{j=1}^p(u_j^{-1}(0))\setminus Z_s}.$\vspace*{2mm}\\
\textit{Step $i+1.$}
Define $C_{i+1}=\overline{\hat{w}_{i+1}(A_{1})}^z.$ Then $C_{i+1}\varsubsetneq\mathbf{C}^{k+t_{i}}$
is irreducible because
$C_{i+1}\subseteq N_i^{-1}(0)$ and $A_{1}$ is irreducible.
Then by Lemma \ref{monter} there are $\delta_{i+1},
q_{t_i+1},\ldots,q_{t_{i+1}}\in\mathbf{C}[\hat{w}_{i+1}],$ where $t_{i+1}=k+t_i-dim C_{i+1},$
such that
$$C_{i+1}\setminus\delta_{i+1}^{-1}(0)=\{\hat{w}_{i+1}\in\mathbf{C}^{k+t_i}
\setminus\delta_{i+1}^{-1}(0):
q_{1}(\hat{w}_{i+1})=\ldots=q_{t_{i+1}}(\hat{w}_{i+1})=0\},$$
and
$\delta_{i+1}|_{C_{i+1}}\neq 0,$ and for every $a\in C_{i+1}\setminus\delta_{i+1}^{-1}(0)$ the map
$(q_{1},\ldots,q_{t_{i+1}}):\mathbf{C}^{k+t_i}\rightarrow\mathbf{C}^{t_{i+1}}$
is a submersion in some
neighborhood of $a$ in $\mathbf{C}^{k+t_i}.$ Moreover,
$\delta_{i+1} I(C_{i+1})\subseteq I(q_1,\ldots,q_{t_{i+1}}).$
Every irreducible component $Z$ of the variety $\overline{\bigcup_{j=1}^p(u_j^{-1}(0))\setminus Z_s}$
satisfies
$\hat{w}_{i+1}(Z)\nsubseteq C_{i+1}$ because $C_{i+1}\subseteq C_i\times\C^{t_i-t_{i-1}}$ and
$\hat{w}_i(Z)\nsubseteq C_i.$ Therefore, in view of Remark \ref{1222}, we may assume that
every such component $Z$ satisfies
$\hat{w}_{i+1}(Z)\nsubseteq\delta_{i+1}^{-1}(0).$

The inclusion $C_{i+1}\subseteq N_i^{-1}(0)$ implies
$\delta_{i+1}N_i=\sum_{j=1}^{t_{i+1}}
q_jr_{i+1,j},$
where $r_{i+1,j}\in\mathbf{C}[\hat{w}_{i+1}],$
and the fact that $u_1$ is a minimal defining function implies that there is
$v_{j}\in\mathcal{O}({D})$ such that
$q_{j}(\hat{w}_{i+1}(x))=v_{j}(x)u_1(x)$
for $j=t_i+1,\ldots,t_{i+1}.$

Let $\hat{v}_{i+1}$
denote the tuple $(v_{t_i+1},\ldots,v_{t_{i+1}})$
of $t_{i+1}-t_i$ variables. Define
$N_{i+1}\in\mathbf{C}[\hat{w}_{i+1},\hat{v}_{i+1}]$ by the formula
$$N_{i+1}(\hat{w}_{i+1},\hat{v}_{i+1})=\sum_{j=1}^{t_{i+1}}v_{j}
r_{{i+1},j}(\hat{w}_{i+1})$$
and observe that
$$N_{i+1}(\hat{w}_{i+1}(x),\hat{v}_{i+1}(x))
=u_1(x)^{k_1-i-1}u_2(x)^{k_{2,i+2}}\cdot\ldots\cdot
u_p(x)^{k_{p,i+2}}R_{i+1}(x),$$
where $R_{i+1}\in\mathcal{O}(D)$
satisfies
$A_{l}\nsubseteq R_{i+1}^{-1}(0)$ for $l=1,\ldots,p.$

Let
$V_{i+1}\subset\mathbf{C}^{k+t_{i+1}+p+1}$
be the algebraic variety defined by the system of equations (in the variables
$\hat{w}_{i+1},\hat{v}_{i+1},u_1,\ldots,u_p, R_{i+1}$):\vspace*{2mm}\\
(E,$i+1$)\hspace*{14.5mm}$N_{i+1}(\hat{w}_{i+1},\hat{v}_{i+1})=
u_1^{k_{1}-i-1}u_2^{k_{2,i+2}}\ldots u_p^{k_{p,i+2}}R_{i+1},$\\
(F,$j$)\hspace*{21mm}$q_{j}(\hat{w}_{i+1})
=v_{j}u_1, \mbox{ }$ for $j=1,\ldots,t_{i+1}.$ \vspace*{2mm}\\
(For $j=1,\ldots,t_i,$ the polynomial $q_j$ was defined in
previous steps.)
Put $\hat{w}_{i+2}(x)=(\hat{w}_{i+1}(x),\hat{v}_{i+1}(x))$
and define $g_{i+1}\in\mathcal{O}(D)$ by
$$g_{i+1}(x)=(\hat{w}_{i+2}(x),u_1(x),\ldots,u_p(x),R_{i+1}(x)).$$
If $k_{1}=i+1$
then $\tilde{f}_{D}=g_{i+1}, \tilde{V}=\overline{g_{i+1}(E)}^z\subseteq V_{i+1}$ satisfy the
requirements
(see Claims \ref{property1}, \ref{property2}). Otherwise we go to Step i+2.\vspace*{0mm}
\begin{claim}\label{property1}The following hold:\\
$(1)$ $\overline{{g}_{k_1}(E)}^z$ is an irreducible component of $V_{k_1}.$
In particular, $$\mathrm{Sing}(\overline{g_{k_1}(E)}^z)
\subseteq\mathrm{Sing}(V_{k_1}).$$
$(2)$ Every $(n-1)$--dimensional irreducible component $S$ of $g_{k_1}^{-1}(\mathrm{Sing}(V_{k_1}))$
with $S\cap E\neq\emptyset,$ satisfies $$S\subseteq\overline{G_D^{-1}(\mathrm{Sing}(N^{-1}(0)))
\setminus A_1}.$$ In particular, after shrinking $D$ if necessary, we have
$${g}_{k_1}^{-1}(\mathrm{Sing}(V_{k_1}))_{(n-1)}
\subseteq\overline{G_{D}^{-1}(\mathrm{Sing}(N^{-1}(0)))\setminus A_{1}}.$$
\end{claim}
\proof Recall that $Z_s=G_{D}^{-1}(\mathrm{Sing}(N^{-1}(0)))_{(n-1)}$ and
$Z_r=\overline{G_{D}^{-1}(N^{-1}(0))\setminus Z_s}$ and suppose that
there is an $(n-1)$--dimensional
irreducible component $S$ of $g_{k_1}^{-1}(\mathrm{Sing}(V_{k_1}))$
with $S\cap E\neq\emptyset$ such that
$S\nsubseteq\overline{Z_s\setminus A_{1}}.$
We consider two
cases:\vspace{2mm}\\
(a) $S\subseteq\bigcup_{j=2}^pu^{-1}_j(0),$\\
(b) $S\nsubseteq\bigcup_{j=2}^pu^{-1}_j(0),$\vspace*{2mm}\\
to show that there is $a\in S$ such that $g_{k_1}(a)\in\mathrm{Reg}(V_{k_1}),$
which contradicts the hypothesis.

Let us begin with (a). The properties of $u_1(x),\ldots,u_p(x)$ imply that for the generic
$a\in S,$ $u_1(a)\neq 0$ hence (in a neighborhood of $g_{k_1}(a)$) the system (F,$j$),
$j=1,\ldots,t_{k_1},$ depicts the graph of the rational map
$(u_1,\hat{w}_1)
\mapsto({v}_{1}(u_1,\hat{w}_1),\ldots,
{v}_{t_{k_1}}(u_1,\hat{w}_1)).$ Moreover, the definition of $N_i$
and the equations (E,$k_1$), (F,$1$),..., (F,$t_{k_1}$) and
$\delta_{i}N_{i-1}=\sum_{j=1}^{t_i}q_jr_{i,j},$ for $i=1,\ldots,k_1$ (where $N_0=N$), imply
that for the generic $a\in S,$ (in a neighborhood of $g_{k_1}(a)$) the variety $V_{k_1}$ is
described by: (F,$j$), $j=1,\ldots,t_{k_1},$\vspace*{2mm}\\
(z)\hspace*{15mm}$N(\hat{w}_1)\prod_{i=1}^{k_1}\delta_i(\hat{w}_i)=
u_1^{k_1}u_2^{k_{2,k_1+1}}\cdot\ldots\cdot u_p^{k_{p,k_1+1}}R_{k_1}.$\vspace*{2mm}\\
Now using (F,$j$)
we can eliminate
the variables $v_1,\ldots,v_{t_{k_1}}$ from (z)
to obtain\vspace*{2mm}\\
(*)\hspace*{20mm}$F(u_1,\hat{w}_1)N(\hat{w}_1)=u_2^{k_{2,k_1+1}}
\ldots u_p^{k_{p,k_1+1}}R_{k_1},$\vspace*{2mm}\\
where $F$ is rational, $(u_1(a),G_{D}(a))\in\mathrm{dom}F,$
and $F(u_1(a),G_{D}(a))\neq 0$ for the gene\-ric $a\in S.$ (The last property due to
$S\subseteq\overline{\bigcup_{j=1}^pu_j^{-1}(0)\setminus Z_s}$ which implies
$\hat{w}_i(S)\nsubseteq\delta_i^{-1}(0)$ for $i=1,\ldots,k_1.$)

Let $\hat{V}$ denote the set defined by (*). To complete
the case (a) it is sufficient to show that $\hat{g}(a)\in\mathrm{Reg}{(\hat{V})},$
for the
generic $a\in S,$
where $\hat{g}$ is the map consisting of those components of $g_{k_1}$
which correspond to the variables appearing in (*).
By assumptions $N(G_{D}(a))=0$ for every $a\in S.$ Moreover,
by the facts that $S\nsubseteq Z_s$ and
$N$ is reduced, there is
$j\in\{1,\ldots,k\}$ such that $\frac{\partial N}{\partial w_j}|_{G_D(S)}\neq 0$ which
clearly implies that $\hat{g}(a)\in\mathrm{Reg}(\hat{V})$ for the generic $a\in S.$

Let us turn to (b). Let $\tilde{g}$
be the map consisting of those components of
$g_{k_1}$ which correspond to the variables appearing in (E,$k_1$)
and let $\bar{g}$
be the map consisting of those components of
$g_{k_1}$ which correspond to the variables appearing in (F,$j$)
for $j=1,\ldots,t_{k_1}.$
For the generic $a\in S,$
the equation  (E,$k_1$) depicts, in a neighborhood of $\tilde{g}(a),$
the graph of the rational function
$R_{k_1}=R_{k_1}(\hat{w}_{k_1},\hat{v}_{k_1},u_2,\ldots,u_p),$
and the variable
$R_{k_1}$ does not appear in any of (F,$j$), $j=1,\ldots,t_{k_1}.$
Hence if we show that for every $a$ in an open dense subset of $S,$
the system of equations (F,$j$), for $j=1,\ldots,t_{k_1},$ defines a manifold in
a neighborhood of $\bar{g}(a)$, then we obtain a contradiction with the assumption
that $g_{k_1}(S)\subseteq\mathrm{Sing}(V_{k_1}).$

We have two cases. If $S\nsubseteq u_1^{-1}(0),$ there is nothing to
prove because each of the considered equations can be divided by $u_1.$
If $S\subseteq u_1^{-1}(0),$ then
for the generic $a\in S,$ the map
$(u_1,\hat{w}_{k_1},\hat{v}_{k_1})\mapsto(q_1(\hat{w}_{k_1})-v_1u_1,
\ldots,q_{t_{k_1}}(\hat{w}_{k_1})-v_{t_{k_{1}}}u_1)$
is a submersion in a neighborhood of $\bar{g}(a).$ This is because
$(q_1,\ldots,q_{t_{k_1}})$ is a submersion in a neighborhood
of $\hat{w}_{k_1}(a),$ and $u_1(a)=q_j(\hat{w}_{k_1}(a))=0,$ for $j=1,\ldots,t_{k_1}.$

It remains to check that $\overline{g_{k_1}(E)}^z$ is an irreducible component of $V_{k_1}.$
We know that $\mathrm{dim}(\overline{g_{k_1}(E)}^z)\geq k+p$ because $u_1,\ldots,u_p$ has been
chosen in such a way that $\overline{(G,u_1,\ldots,u_p)(E)}^z=\C^{k+p}.$ So it is sufficient to
check that there is $a\in E$ such that $V_{k_1}$ is a $(k+p)$--dimensional
manifold in some neighborhood of $g_{k_1}(a)$. But this holds for $a\in E$ with
$u_1(a)\cdot\ldots\cdot u_p(a)\neq 0.$ Indeed, then in a neighborhood of $g_{k_1}(a),$
$V_{k_1}$ is the graph of the rational map $(u_1,\ldots,u_p,w_1,\ldots,w_k)\mapsto(v_1,
\ldots,v_{t_{k_1}},R_{k_1}).$\qed
\begin{claim}\label{property2} If there is a sequence
$g_{k_1,\nu}\in\mathcal{O}(E,V_{k_1})$ of Nash maps
converging uniformly to $g_{k_1}|_E$ then there are a sequence
$G_{\nu}\in\mathcal{O}({E},\mathbf{C}^{k})$
of Nash maps converging uniformly to $G$ and an open neighborhood $D'$ of $E$
such that $\{(N\circ G_{\nu,D'})^{-1}(0)\}$ converges to $(N\circ G_{D'})^{-1}(0)$
in the sense of chains.
\end{claim}
\proof Let $g_{k_1,\nu, D'}\in\mathcal{O}(D',V_{k_1})$ be a sequence of Nash maps converging
uniformly to $g_{k_1,D'},$ where $D'$ is an open neighborhood of $E$ in $D$ such that
$(N\circ G_{D'})^{-1}(0)=(A_1\cup\ldots\cup A_p)\cap D',$ for $A_1,\ldots,A_p$
introduced in the proof of Proposition \ref{hhhk}.

The map
$g_{k_1,\nu,D'}$ is of the form:
$$g_{k_1,\nu,D'}(x)=(\hat{w}_{k_1+1,\nu}(x),
u_{1,\nu}(x),\ldots,u_{p,\nu}(x),R_{k_1,\nu}(x)),$$
where $\hat{w}_{i+1,\nu}(x)=(\hat{w}_{i,\nu}(x),\hat{v}_{i,\nu}(x)),$
$\hat{v}_{i,\nu}(x)=(v_{t_{i-1}+1,\nu}(x),\ldots,v_{t_i,\nu}(x)),$
for $i=1,\ldots,k_1,$ where $t_0=0.$
We check
that $G_{\nu,D'}(x)=\hat{w}_{1,\nu}(x)$ satisfies the requirements.

Clearly, it is sufficient to show that for every $l\in\{1,\ldots,p\}$
and for
the generic point $a\in A_l\cap D'$ there is a neighborhood $U$ of $a$
in $D'$ such that $\{(N\circ G_{\nu,D'})^{-1}(0)\cap U\}$ converges to
$A_l\cap U$ in the sense of chains. Fix $l\in\{1,\ldots,p\}.$

The components of $g_{k_1,\nu,D'}$ satisfy the equation (E,$k_1$)
therefore\vspace*{2mm}\\
(E,$k_1,\nu$)\hspace*{5mm}$
N_{k_1}(\hat{w}_{k_1+1,\nu}(x))
=u_{2,\nu}(x)^{k_{2,k_1+1}}\ldots u_{p,\nu}(x)^{k_{p,k_1+1}}R_{k_1,\nu}(x),$\vspace*{2mm}\\
for every $x\in D', \nu\in\mathbf{N}.$

By the definition of $N_{i+1}$ and by the fact that the components of $g_{k_1,\nu,D'}$
satisfy the equations $\delta_{i+1}N_i=\sum_{j=1}^{t_{i+1}}q_jr_{i+1,j},$
$(F,1),\ldots, (F,t_{k_1}),$ we have
$$N_{i+1}(\hat{w}_{i+2,\nu}(x))u_{1,\nu}(x)=\delta_{i+1}(\hat{w}_{i+1,\nu}(x))
N_i(\hat{w}_{i+1,\nu}(x)),$$
for $i=0,\ldots,k_1-1,$ $x\in D',$ where $N_0(\hat{w}_{1,\nu}(x))=N(\hat{w}_{1,\nu}(x)).$
This implies that
$$N_{k_1}(\hat{w}_{k_1+1,\nu}(x))u_{1,\nu}(x)^{k_1}=
\tilde{T}_{\nu}(x)N(\hat{w}_{1,\nu}(x)),$$
for some $\tilde{T}_{\nu}\in\mathcal{O}(D'),$ which combined with (E,$k_1,\nu$)
gives\vspace*{2mm}\\
($\alpha$)\hspace*{10mm}$\tilde{T}_{\nu}(x)N(\hat{w}_{1,\nu}(x))=
u_{1,\nu}(x)^{k_1}u_{2,\nu}(x)^{k_{2,k_1+1}}\ldots u_{p,\nu}(x)^{k_{p,k_1+1}}R_{k_1,\nu}(x),$\vspace*{2mm}\\
for every $x\in D',\nu\in\mathbf{N}.$

The facts that
$A_{l}\subseteq\overline{u_l^{-1}(0)\setminus R_{k_1}^{-1}(0)}$
and that
$\mathrm{dim}(u_l^{-1}(0)\cap u_t^{-1}(0))<n-1,$ for every $t\neq l,$
clearly imply that for the generic $a\in A_l\cap D'$
there is an open neighborhood $U\Subset D'$ such that
$(\bigcup_{j\neq l}u_j^{-1}(0)\cup R_{k_1}^{-1}(0))\cap\overline{U}=\emptyset.$ Consequently,
for sufficiently
large $\nu,$ by ($\alpha$), we have
$$(\tilde{T}_{\nu}(x)N(\hat{w}_{1,\nu}(x)))^{-1}(0)\cap U=
u_{l,\nu}^{-1}(0)\cap U.$$
Now by the fact that $u_l$ is a minimal defining function,
$\{u_{l,\nu}^{-1}(0)\cap U\}$ converges to $A_l\cap U$ in the sense of chains.
On the other hand, $\{(N(\hat{w}_{1,\nu}(x)))^{-1}(0)\cap U\}$ converges to $A_l\cap U$
locally uniformly and, in view of the last equation, for $\nu$ large enough,
$(N(\hat{w}_{1,\nu}(x)))^{-1}(0)\cap U\subseteq u_{l,\nu}^{-1}(0)\cap U$ so
$\{(N(\hat{w}_{1,\nu}(x)))^{-1}(0)\cap U\}$ converges
to $A_l\cap U$ in the sense of chains.\qed\vspace*{2mm}

Once we have proved Claims \ref{property1}, \ref{property2},
the proof of Proposition \ref{hhhk} is also completed.\qed

\section{Generalization of Theorem \ref{maintheorem}}\label{general21}
Let $f:U\rightarrow V$ be as in Theorem \ref{maintheorem}. As already mentioned, without loss
of generality, we can additionally assume in Theorem \ref{maintheorem} that
$\overline{f(U)}^z=V$ and $V\subset\C^m\times\C^k$ has proper projection onto $\C^m$
where $m=\mathrm{dim}(V)$. {Write $f=(\tilde f, \hat f)$, where}
$\tilde{f}, \hat{f}$ denote the first $m$ and the last $k$ components of $f,$ respectively. Then
$\tilde{f}(U)\nsubseteq\Sigma_V.$ 

Let $\mathcal V(\tilde f)$ 
 denote the pull-back of $V$ by $\tilde f,$ i. e. $\mathcal V(\tilde f)=(\tilde{f}\times\mathrm{id_{\C^k}})^{-1}(V)$. Then the fact that $f(U)\subseteq V$ can be equivalently stated as 
$\mathrm{graph}(\hat{f})\subset\mathcal V(\tilde f).$
Under these assumptions Theorem \ref{maintheorem} can be reformulated as
follows. \vspace*{2mm}\\
\textbf{Theorem 2.1'} \emph{For every open $\tilde{U}\Subset U$
there are a sequence
$\tilde{f}_{\nu}:\tilde{U}\rightarrow\C^m$ of Nash maps converging uniformly to
$\tilde{f}|_{\tilde{U}}$ and a sequence $M_{\nu}$ of Nash sets of pure dimension $n$
converging to
$\mathrm{graph}(\hat{f})\cap (\tilde{U}\times\C^k)$ in the sense of chains such that
$M_{\nu}\subset
\mathcal V(\tilde f_\nu)$
for every $\nu$.}\vspace*{2mm}\\
 {Indeed, in this case, for $\nu$ large, (shrinking $\tilde{U}$ slightly we obtain that) $M_\nu$ is the graph of a map 
that defines the second 
 part $\hat f_\nu$ of $f_\nu = (\tilde f_\nu, \hat f_\nu)$, cf. the proof of Corollary \ref{maintool}.  
 But the method of the proof gives that}    
 $\tilde{f}$ can be approximated by Nash maps $\tilde{f_{\nu}}$ in such a way 
that all purely $n$--dimensional analytic sets (in particular all graphs of maps
holomorphic on $U$) 
contained in $\mathcal V(\tilde f)$ can be simultaneously
approximated in the sense of chains by Nash sets contained in 
$\mathcal V(\tilde f_\nu)$.
More precisely, the following generalization of Theorem~2.1' holds.
\begin{theorem}\label{generalmain}
 Let $U$ be an open polydisc in $\mathbf{C}^n$ and let $V\subset\mathbf{C}^m\times\mathbf{C}^k$
be an algebraic variety of pure dimension $m$ with proper projection onto $\mathbf{C}^m.$ 
Let $\tilde f :U\rightarrow\C^m$ be a holomorphic map such that $\tilde{f}(U)\nsubseteq\Sigma_V.$ Then  
for every open $\tilde{U}\Subset U$
there is a sequence
$\tilde{f}_{\nu}:\tilde{U}\rightarrow\C^m$ of Nash maps converging uniformly to
$\tilde{f}|_{\tilde{U}}$ such that for every analytic set $M\subset
\mathcal V(\tilde f)$
of pure dimension $n$ there is a sequence $M_{\nu}$ of Nash sets of pure dimension $n$
converging to
$M\cap (\tilde{U}\times\C^k)$ in the sense of chains such that
$M_{\nu}\subset \mathcal V(\tilde f_\nu)$ for every $\nu.$
\end{theorem}
\proof Let $N_1,\ldots, N_s$ be polynomials in $m$ 
complex variables such that $\Sigma_V=\{N_1=\ldots=N_s=0\}.$
Shrinking $U$ if needed we can assume that $\tilde{f}^{-1}(\Sigma_V)$ has finitely many, say $p,$
$(n-1)$-dimensional irreducible components. Denote these components by $C_1,\ldots,C_p.$ 
Let $u_1,\ldots,u_p$ be minimal defining
functions for $C_1,\ldots,C_p,$ respectively. Then there are $R_j\in\mathcal{O}(U)$ and $k_{j,i}\in\mathbf{N},$
for $j=1,\ldots,s$ and $i=1,\ldots,p,$ 
such that $N_j\circ\tilde{f}=R_ju_{1}^{k_{j,1}}\cdot\ldots\cdot u_p^{k_{j,p}},$ 
and every $R_j$ does not vanish identically on any $C_i.$ 

Fix open $\tilde{U}, \hat{U}$ with $\tilde{U}\Subset\hat{U}\Subset U.$
By Theorem \ref{maintheorem}, there are Nash maps
$\tilde{f}_{\nu}, R_{j,\nu}, u_{i,\nu}$ approximating $\tilde{f}, R_j, u_{i},$ respectively, on
$\hat{U}$ and such that $N_j\circ\tilde{f}_{\nu}=R_{j,\nu}u_{1,\nu}^{k_{j,1}}\cdot\ldots\cdot u_{p,\nu}^{k_{j,p}},$ for $j=1,\ldots,s.$ Now $X=\mathcal{V}(\tilde{f})\cap(\hat{U}\times\C^k),$
$X_{\nu}=\mathcal{V}(\tilde{f}_{\nu})$ satisfy the assumptions of Theorem \ref{component}.
Application of this theorem completes the proof.\qed\vspace*{2mm}\\
{\small\textit{Acknowledgements.} We wish to thank Professor Wojciech Kucharz for
reading an early version of this manuscript and for helpful comments. We also thank Dr Jakub
Byszewski for a helpful discussion.}

{\small
\begin{thebibliography}{}
\bibitem{An}
Andr\'e, M.: \emph{Cinq expos\'es sur la d\'esingularization,}
manuscript, \'Ecole Polytechnique F\'ed\'erale de Lausanne, 1992
\bibitem{A68}
Artin, M.: \emph{On the solutions of analytic equations,} Invent. Math.  \textbf{5}, 277--291 (1968)
\bibitem{A}
Artin, M.: \emph{Algebraic approximation of structures over complete
local rings,} Publ. I.H.E.S. \textbf{36}, 23-58 (1969)
\bibitem{Ba}
Barlet, D.: \emph{Espace analytique r\'eduit des cycles
analytiques complexes compacts d'un espace analytique complexe de
dimension finie,} Fonctions de plusieurs variables complexes, II,
S\'em. Fran\c{c}ois Norguet, 1974-1975, Lecture Notes in Math.,
\textbf{482}, pp. 1-158, Springer, Berlin 1975
\bibitem{B2}
Bilski, M.:
\emph{Approximation of analytic sets with proper projection by algebraic sets},
Constr. Approx. \textbf{35}, 273--291 (2012).
\bibitem{B4}
Bilski, M.: \emph{Algebraic approximation of analytic sets and
mappings,} J. Math. Pures Appl. \textbf{90}, 312-327 (2008)
\bibitem{BoK2}
Bochnak, J., Kucharz, W.: \emph{Approximation of holomorphic maps
by algebraic morphisms,} Ann. Polon. Math. \textbf{80}, 85-92
(2003)
\bibitem{Ch}
Chirka, E.: Complex analytic sets. Kluwer Academic Publ.,
Dordrecht-Boston-London 1989
\bibitem{DLS}
Demailly, J.-P., Lempert, L., Shiffman, B.: \emph{Algebraic
approximation of holomorphic maps from Stein domains to projective
manifolds,} Duke Math. J. \textbf{76}, 333-363 (1994)
\bibitem{vD}
van den Dries, L.: \emph{A specialization theorem for analytic functions
on compact sets,} Nederl. Akad. Wetensch. Indag. Math. \textbf{44},
391-396 (1982)
\bibitem{DF1}
Drinovec Drnov\v sek, B., Forstneri\v c, F.: \emph{Strongly pseudoconvex domains
as subvarieties of complex manifolds.} Amer. J. Math. \textbf{132}, 331-360
(2010)
\bibitem{DF2}
Drinovec Drnov\v sek, B., Forstneri\v c, F.: \emph{Holomorphic curves in complex spaces.}
Duke Math. J. \textbf{139}, 203-253 (2007)
\bibitem{FT}
Fatabbi, G., Tancredi, A.: \emph{On the factoriality of some rings of
complex Nash functions.} Bull. Sci. Math. \textbf{126}, 61-70 (2002)
\bibitem{Fo}
Forstneri\v c, F.: \emph{Holomorphic flexibility properties of
complex manifolds,} Amer. J. Math. \textbf{128}, 239-270 (2006)
\bibitem{Fr}
Frisch, J.: \emph{Points de platitude d'un morphisme d'espaces analytiques
complexes.} Invent. Math. \textbf{4}, 118-138 (1967)
\bibitem{Hartshorne77}
Hartshorne, R.: Algebraic geometry. Graduate Texts in Mathematics, No. 52. Springer-Verlag, New York, 1977 
\bibitem{HR}
Hauser, H., Rond, G.: \emph{Artin Approximation.} Preprint:\\ 
homepage.univie.ac.at/herwig.hauser/
\bibitem{Ho}
H\"{o}rmander, L.: An introduction to complex analysis in several
variables. North-Holland, Amsterdam-New York-Oxford-Tokyo 1990
\bibitem{Ku}
Kucharz, W.: \emph{The Runge approximation problem for holomorphic
maps into Grassmannians,} Math. Z. \textbf{218}, 343-348 (1995)
\bibitem{Lan}
Lang, S.: Algebra. Revised third edition. Graduate Texts in Mathematics, 211. Springer-Verlag,
New York, 2002
\bibitem{Lem}
Lempert, L.: \emph{Algebraic approximations in analytic geometry,}
Invent. Math. \textbf{121}, 335-354 (1995)
\bibitem{Lo}
\L ojasiewicz, S.: Introduction to complex analytic geometry.
Birkh\"{a}user, Basel 1991
\bibitem{MY}
Marinescu, G., Yeganefar, N.: \emph{Embeddability of some strongly pseudoconvex
CR manifolds,} Trans. Amer. Math. Soc. \textbf{359}, 4757-4771 (2007)
\bibitem{Mi}
Mir, N.: \emph{Algebraic approximation in CR geometry.}
J. Math. Pures Appl. \textbf{98}, 72-88 (2012)
\bibitem{NR}
Napier, T., Ramachandran, M.: \emph{The $L^2 \bar{\partial}$--method, weak Lefschetz
theorems, and the topology of K\"ahler manifolds,} J. Amer. Math. Soc. \textbf{11},
375-396 (1998)
\bibitem{Og}
Ogoma, T.: \emph{General N\'eron desingularization based on the
idea of Popescu.} J. of Algebra \textbf{167}, 57-84 (1994)
\bibitem{Po1}
Popescu, D.: \emph{General N\'eron desingularization.} Nagoya
Math. J. \textbf{100}, 97-126 (1985)
\bibitem{Po2}
Popescu, D.: \emph{General N\'eron desingularization and
approximation.} Nagoya Math. J. \textbf{104}, 85-115 (1986)
\bibitem{Sp}
Spivakovsky, M.: \emph{A new proof of D. Popescu's theorem on
smoothing of ring homomorphisms.} J. Amer. Math. Soc.,
\textbf{12}, 381-444 (1999)
\bibitem{Sn}
Swan, R.: N\'eron-Popescu desingularization. Algebra and geometry (Taipei, 1995), 135-192,
Lect. Algebra Geom., 2, Internat. Press, Cambridge, MA, 1998 
\bibitem{TT3}
Tancredi, A., Tognoli, A.: \emph{On the relative Nash approximation of
analytic maps,} Rev. Mat. Complut. \textbf{11}, 185-201 (1998)
\bibitem{Tw}
Tworzewski, P.: \emph{Intersections of analytic sets with linear
subspaces,} Ann. Sc. Norm. Super. Pisa \textbf{17,} 227-271 (1990)
\bibitem{Tw2}
Tworzewski, P.: \emph{Intersection theory in complex analytic geometry,}
Ann. Polon. Math., \textbf{62.2} 177-191 (1995)
\bibitem{TwW}
Tworzewski, P., Winiarski, T.: \emph{Continuity of intersection of
analytic sets,} Ann. Polon. Math. \textbf{42,} 387-393 (1983)
\end{thebibliography}
}
\end{document}